\documentclass[hidelinks]{article}
\usepackage[utf8]{inputenc}
\usepackage[utf8]{inputenc}
\usepackage{mathrsfs}
\usepackage{amssymb}
\usepackage[margin=1in]{geometry}
\usepackage{changepage}
\usepackage{stmaryrd}
\usepackage{xcolor}
\usepackage{bbold}
\usepackage{amsmath,amsthm}
\usepackage{graphicx}
\usepackage{hyperref}
\usepackage{tikz}
\usepackage{soul}
\newtheorem{propo}{Proposition}
\newtheorem{remark}{Remark}
\newtheorem{lemma}{Lemma}

\newtheorem{theorem}{Theorem}
\newtheorem*{theorem*}{Theorem}

\numberwithin{equation}{section}
\newcommand{\vertiii}[1]{{\left\vert\kern-0.25ex\left\vert\kern-0.25ex\left\vert #1 
		\right\vert\kern-0.25ex\right\vert\kern-0.25ex\right\vert}}

\def\N{\mathbb{N}}
\def\R{\mathbb{R}}

\def\d{\,\mathrm{d}}

\def\E{\mathbb{E}}

\def\eps{\varepsilon}

\title{\textbf{Weak Error on the densities for the Euler scheme of stable additive SDEs with  H\"older drift}}
\author{Mathis Fitoussi\footnote{Laboratoire de Mathématiques et Modélisation d'Evry (LaMME), UMR CNRS 8071, Université d'Evry Val d'Essonne-Paris Saclay, 23 Boulevard de France, 91037 Evry, mathis dot fitoussi at univ-evry dot fr}, Stéphane Menozzi\footnote{Laboratoire de Mathématiques et Modélisation d'Evry (LaMME), UMR CNRS 8071, Université d'Evry Val d'Essonne-Paris Saclay,  23 Boulevard de France, 91037 Evry, stephane dot menozzi at univ-evry dot fr}}
\date{\today}

\begin{document}
\maketitle


\begin{center}
	\textbf{Abstract}
\end{center}
\begin{adjustwidth}{0.8in}{0.8in}
    For a fixed finite time horizon $T>0$, we are interested in  the Euler-Maruyama dicretization of the SDE
\begin{equation}\label{sde-abstract}\notag
\d X_t = b(t,X_t)\d t + \d Z_t, \qquad X_0 = x \in \R^d,
\end{equation}
where $Z_t$ is a symmetric isotropic $d$-dimensional $\alpha$-stable process, $\alpha\in (1,2] $ and the drift $b\in L^\infty \left([0,T],\mathcal{C}^{\beta}(\R^d,\R^d)\right)$, $\beta\in (0,1)$, is bounded and H\"older regular in space. Using an Euler scheme with a randomization of the time variable, we show that, denoting $\gamma:=\alpha+\beta-1$, the weak error on densities related to this discretization converges at the rate $\gamma/\alpha$.
\end{adjustwidth}

\section{Introduction}
For a fixed finite time horizon $T>0$, we are interested in  the Euler-Maruyama dicretization of the SDE
\begin{equation}\label{hold-sde}
	\d X_t = b(t,X_t)\d t + \d Z_t, \qquad X_0 = x, \qquad \forall t \in [0,T],
\end{equation}
where $Z_t$ is a symmetric isotropic $d$-dimensional $\alpha$-stable process, $\alpha\in (1,2] $ and $b\in L^\infty \left([0,T],\mathcal{C}^{\beta}(\R^d,\R^d)\right)$, $\beta\in (0,1)$, {i.e. it} is bounded and H\"older regular in space. {In this setting, weak well-posedness holds for \eqref{hold-sde} since the \textit{natural} condition 
	\begin{equation}\label{hold-serrin}
		\gamma:=\beta+\alpha-1>0 \iff \beta + \alpha > 1,
	\end{equation}
	is  always satisfied. The condition \eqref{hold-serrin} actually ensures weak well-posedness for the SDE \eqref{hold-sde}, even in the super-critical case $\alpha\in (0,1]$, provided the drift is time homogeneous or bounded in time (see \cite{TTW74}, \cite{MP14}, \cite{CZZ21}, see also \cite{Pri12}, \cite{CZZ21} for strong well-posedness established under the more stringent condition $\beta + \alpha/2 >1$)}.\\

The goal of this chapter is to prove a convergence rate for the weak error on densities associated with an {\textit{appropriate}} Euler scheme for \eqref{hold-sde}.
\subsection{Definition of the scheme}\label{hold-subsec-def-scheme}

We will use a discretization scheme with $n$ time steps over $[0,T]$, with constant step size $h:=T/n$. For the rest of this chapter, we denote, $\forall k \in \{{0},...,n\}, t_k := kh$ and $\forall s>0, \tau_s^h := h \lfloor \frac{s}{h} \rfloor \in (s-h,s]$, which is the last grid point before time $s$. Namely, if $s\in [t_k,t_{k+1}), \tau_s^h = t_k$. \\

We define a step of the Euler scheme, {starting from $X_0=x$}, as
\begin{equation}\label{hold-euler-scheme}
	X_{t_{k+1}}^h = X_{t_k}^h + hb(U_k,X_{t_k}^h)+ (Z_{t_{k+1}}-Z_{t_{k}}),\ k\in \N,
\end{equation}
where the $(U_k)_{k\in \N} $ are independent random variables, independent as well from the driving noise, s.t. $U_k\overset{({\rm law})}={\mathcal U}([{t_k,t_{k+1}}]) $, {i.e. $U_k $ is uniform on the time interval $[t_k,t_{k+1}]$}.
We consider the corresponding  time interpolation defined as the solution to
\begin{equation}\label{hold-scheme-interpo}
	\d X^h_t=b(U_{\tau_t^h/h},X^h_{\tau^h_t}) \d t+ \d Z_t.
\end{equation}
As $b$ is bounded, equation \eqref{hold-scheme-interpo} is well-posed and $X_t^h$ admits a density for $t>0${. We} can refer to this end to \cite{FJM24} for related estimates. We will denote by $\Gamma^h(0,x,t,\cdot)$ this density at time $t\in (0,T]$ when starting from $x$ at time $0$.


\subsection{Euler scheme - state of the art}

For all $0\leq s < t \leq T$, it is known that the unique weak solution to \eqref{hold-sde}, starting in $x$ at time $s$ admits a density, which we will denote $\Gamma (s,x,t,\cdot)$.  It has as well been established in \cite{MZ22} that in the {current} setting $\Gamma $ enjoys two-sided stable heat kernel estimates for $\alpha\in (1,2) $ whereas this property can already be derived from Friedman \cite{frie:64} {(under some additional smoothness in time for $b$) or \cite{MPZ21}} in the Brownian case. In this chapter, we are interested in the weak error on densities, which is defined as the {quantity} 
\begin{equation}\label{hold-DIFF_DENS}
	|\Gamma(s,x,t,y)-\Gamma^h(s,x,t,y)|.
\end{equation}
{In particular we want to bound it, up to a multiplicative constant,  by the product of an appropriate power of the time step $h$ and a density which provides an upper bound for the one of  the driving noise. This would then in particular allow to integrate against possibly irregular test functions having the corresponding convergence rate}.\\

The general definition of the weak error is
\begin{equation}
	\label{hold-def-weak-error}\mathcal{E}(f,t,x,h):= \E_{0,x}\left[f(X_t^h)-f(X_t)\right],
\end{equation}
for $f$ belonging to a suitable class of test functions, and where the meaning of the expectation subscript for the rest of the chapter is $\E_{0,x}[\cdot]:=\E[\cdot|X^h_0=X_0=x]$.\\

Deriving convergence results for the weak error involves studying the PDE
\begin{equation}\label{hold-pde}
	\left( \partial_s + b(s,x)\cdot \nabla_x + \mathcal{L}^\alpha \right)u(s,x)={0} \; \mathrm{on} \;[0,T)\times \R^d, \qquad u(T,\cdot) = {f} \; \mathrm{on} \; \R^d,
\end{equation}
where $\mathcal{L}^\alpha$ is the generator of the noise. 
When the coefficients of \eqref{hold-sde} and the test function $f$ are smooth, the seminal paper of Talay and Tubaro (\cite{TT90}) gives a convergence {rate} of order 1 in $h$ in the Brownian case. {Similar results were obtained for the densities in \cite{KM02} and \cite{KM10} respectively in the Brownian and pure-jump settings}. With $\beta$-H\"older coefficients {and again a smooth $f$ in \eqref{hold-def-weak-error}}, the work of Mikulevicius and Platen (\cite{MP91}) proves a convergence in $h^\frac{\beta}{2}$ in the Brownian case. {This result was extended to densities in \cite{KM17}}. In these works, when applying Itô's formula, authors use the regularity of the drift to treat terms of the form $b(r,X_r^h)-b(U_{\tau_r^h/h},X_{\tau_r^h}^h)$ but do not exploit {the full parabolic bootstrap associated with the PDE \eqref{hold-pde}}.
Note that this approach intrinsically leads to a \textit{strong} convergence order {and that all the previously quoted results are for an SDE with multiplicative noise which is as well $\beta $-H\"older continuous in space. In that setting, we believe that the rate is {sharp}. However, for an additive noise as in \eqref{hold-sde} (or a multiplicative noise with \textit{smooth} diffusion coefficient), the rate can be significantly  improved}. \\

%

In the Brownian setting, one way to {proceed} is to use the stochastic sewing lemma introduced in \cite{Le20}, which allows to quantify the discretization error along rough functionals of the Brownian path. In the specific case of a {$\beta $-}H\"older continuous drift and terminal condition {$f$}, in the work \cite{Hol22}, the author improves the convergence rate from \cite{MP91} to $h^{\frac{\beta+1}{2}-{\varepsilon}},\ {\varepsilon>0}$. We aim to extend this result to the pure-jump setting $\alpha\in (1,2]$ and to a more general class of test functions by working on densities {achieving as well $\varepsilon=0 $}. Let us also mention the work \cite{LL22}, which proves a \textit{strong} (i.e. on trajectories) rate of convergence of order $1/2$ (up to a logarithmic factor) in the Brownian setting for {$L_t^q-L_x^p$ drifts under the Krylov-R\"ockner type condition $d/p+2/q<1 $}. 
{However, the use of stochastic sewing techniques still does not allow to take {full advantage of the parabolic bootstrap associated with the fundamental solution of  \eqref{hold-pde} when the test function is rough, e.g. Dirac masses leading to the weak error on densities}.}\\

{In the current work, we precisely focus on these types of errors of the form} $\mathcal{E}(\delta_y,t,x,h)$ (where $\delta_y$ is the Dirac mass at point $y$). {From Itô's formula, \eqref{hold-def-weak-error} and \eqref{hold-pde}, this formally writes}
\begin{align}\label{hold-errdelt}
	\mathcal{E}(\delta_y,t,x,h)=\E _{0,x}\left[\int_0^t \left( b(r,X_r^h)-b(U_{\tau_r^h/h},X_{\tau_r^h}^h)) \right) \cdot \nabla_z \Gamma (r,z,t,y)|_{z=X_r^h} \d r\right].
\end{align}
{To analyze the corresponding error}, a new idea was introduced in \cite{BJ20}. The drift was therein assumed to be merely measurable and bounded so that no rate could be \textit{a priori} derived from the difference in \eqref{hold-errdelt}. The point then  consists in using the regularity of the solution to \eqref{hold-pde} instead of that of $b$. Namely, writing
\begin{align}
	\E_{0,x} [b(r,X_r^h)\cdot\nabla \Gamma (r,X_r^h,t,y)&-b(r,X_{\tau_r^h}^h)\cdot\nabla \Gamma (r,X_{\tau_r^h}^h,t,y) ]\notag\\&= \int [ \Gamma^h (0,x,r,z) - \Gamma^h (0,x,\tau_r^h,z) ]b(r,z) \cdot\nabla \Gamma (r,z,t,y)\d z\label{hold-bj20}
\end{align}
{one can exploit some additional (or-bootstrapped regularity) of $\Gamma^h$ in its forward time variable.
	Namely, it was proved in the Brownian setting of \cite{BJ20} that for a bounded drift, this regularity was of order $1/2 $, which actually formally corresponds to the exponent $\gamma/\alpha $ (with $\gamma $ defined in \eqref{hold-serrin}) when taking $\beta=0 $. This result still holds in the current setting with $\beta\in (0,1) $ and provides a significant  improvement, corresponding to the expected regularity deriving from the parabolic bootstrap in the forward variable, when compared  to the $\beta$-H\"older regularity of $b$ in space}. {To handle the error from \eqref{hold-errdelt} one would need as well to investigate a space sensitivity of the gradient of the density $\Gamma$ in its backward variable. This could as well be done by exploiting the parabolic bootstrap}. \\

On the other hand, we also have to account for terms involving $b(r,X_r)-b(U_{\tau_r^h/h},X_r)$. One way to achieve the expected convergence rate is to make strong assumptions on the time regularity of $b$: we would need $b(\cdot,z)$ to be $\gamma/\alpha$-H\"older. Importantly, without making any assumption on the time regularity of the drift, those terms can be handled thanks to the randomization of the time argument introduced in \eqref{hold-euler-scheme}, which allows for a convenient use of the Fubini theorem in the error analysis (see \eqref{hold-convenient-fubini} below). {This averaging procedure can somehow be seen as well as a regularization by noise phenomenon.}\\


{Let us mention that for the proofs below we will not rely on the previous expansion of the error, which we  presented here in order to give an idea of the main crucial steps and tools for the error analysis, but on the Duhamel representations of the densities expanded using the density of the driving noise as proxy (see Proposition \ref{hold-prop-main-estimates-D} below)}.\\

{From} the above techniques (forward time regularity of $\Gamma^h$ and time randomization), a rate of order $\frac{\alpha - 1 -\left( \frac{d}{p}+\frac{\alpha}{q}\right)}{\alpha}>0$ is derived in \cite{JM211} and \cite{FJM24}, respectively in the Brownian and  pure-jump settings, for a Lebesgue drift in $L^q_t-L^p_x$ {for the difference of the densities \eqref{hold-DIFF_DENS}}. {Comparing this rate to that of \cite{LL22}, although $1/\alpha$ is lost due to the gradient in \eqref{hold-bj20} (time singularity induced by the gradient of the density of the noise), one sees that the convergence rate {displays explicitly 
		the ``gap to singularity" $\alpha-1-(d/p+\alpha/q) $ or Serrin condition in that setting (critical {stable} parabolic scaling in Lebesgue spaces)}.}\\

In Theorem \ref{hold-thm-main}, we derive a  weak error rate in $h^{\frac{\gamma}{\alpha}}$, where $\gamma:= \beta+ \alpha -1$ is the corresponding ``gap to singularity" in the H\"older case. Importantly, if we interpret $-\left( \frac{d}{p}+\frac{\alpha}{q}\right)$ as the regularity in the former works\footnote{actually this exponent naturally appears as the negative regularity parameter when embedding the time-space Lebesgue space in a Besov space with infinite integrability indexes (which can be identified with a usual H\"older space when the regularity index is positive), see e.g. \cite{Sawano18}.}, there is continuity of the rate of the convergence w.r.t. the regularity of the drift. Continuity w.r.t. the stability index $\alpha$ also holds when comparing Theorem \ref{hold-thm-main} to the results in \cite{Hol22} {(and getting rid of the $\varepsilon $ in the rate therein)}, thus extending the former to a more general class of test functions and noises.\\

The restriction to the sub-critical case $\alpha\in (1,2) $, for which the intensity of the driving noise somehow dominates the drift in \textit{small time}, is here required  mainly for approximation purposes. {Again}, in the  \cite{CZZ21} paper, weak well-posedness is obtained {under \eqref{hold-serrin}} without the restriction $\alpha\in (1,2) $. On the other hand, as we are interested in heat kernel estimates, it is also well known, see e.g. Kulik \textit{et al.} \cite{KK18}, \cite{Kul19}, \cite{MZ22}, that in the super-critical regime $\alpha\in (0,1) $ for which the Hölder regularity needs to be large enough to compensate the lower regularizing effects of the noise, considerations on {some flows related to the drift in \eqref{hold-sde}} are needed. These aspects become a rather difficult issue when considering associated discretization schemes (see Konakov \textit{et al.} \cite{KM23} in connection with stochastic algorithms of Robbins Monro type).\\

{The chapter is organized as follows: in Section \ref{hold-subsec-noise} we specify some properties of the driving noise in \eqref{hold-sde}.  
	Section \ref{hold-subsec-main-result} is then dedicated to the statement of the main results (we give some controls on the densities of the SDE and the Euler scheme in Proposition \ref{hold-prop-main-estimates} and the convergence rate for the weak error in Theorem \ref{hold-thm-main}). Section \ref{hold-SEC_PROOF_THM} is devoted to the proof of the main theorem. The proof is achieved via exploiting some additional quantitative properties of the density of the driving noise, the Duhamel representation of the densities (see Proposition \ref{hold-prop-main-estimates-D}) and the regularity results of Proposition \ref{hold-prop-main-estimates} which are in turn proved in Section \ref{hold-SEC_PROOF_REG}.}

\subsection{Driving noise and related density properties}\label{hold-subsec-noise}
Let us denote by $\mathcal{L}^\alpha$ the generator of the driving noise $Z$ and $p_\alpha:\R_+\backslash\{0\}\times \R^d\rightarrow \R_+$ its density. In the case $\alpha = 2$, $\mathcal{L}^\alpha$ is the usual normalized Laplacian $\frac{1}{2}\Delta$. The noise is a Brownian Motion and its gaussian marginal densities are explicit.\\

When $\alpha \in (1,2)$, in whole generality, the generator of a symmetric stable process writes, $\forall \phi \in C_0^\infty (\R^d,\R)$ (smooth compactly supported functions),
\begin{align*}
	\mathcal{L}^\alpha \phi (x)&= \mathrm{p.v.} \int_{\R^d} \left[ \phi(x+z) - \phi(x)\right]\nu(\d z)\\
	&=\mathrm{p.v.}\int_{\R_+}\int_{\mathbb{S}^{d-1}}\left[ \phi(x+\rho \xi) - \phi(x)\right]\mu(\d \xi)\frac{\d \rho}{\rho^{1+\alpha}}
\end{align*}
(see \cite{Sat99} for the polar decomposition of the stable L\'evy measure) where $\mu$ is a symmetric measure on the unit sphere $\mathbb{S}^{d-1}$. We will here restrict to the case where $\mu  =m$  the Lebesgue measure on the sphere but it is very likely that the analysis below can be extended to the case where
$\mu$ is symmetric and $\exists \kappa \geq 1 : \forall \lambda \in \R^d$,
\begin{equation*}
	C^{-1} m(\d \xi) \leq \mu (\d \xi) \leq C m(\d \xi),
\end{equation*} 
i.e. it is equivalent to the Lebesgue measure on the sphere. Indeed, in that setting
Watanabe (see \cite{Wa07}, Theorem 1.5) and Kolokoltsov (\cite{Kol00}, Propositions 2.1--2.5) showed that if $C^{-1} m(\d \xi) \leq \mu (\d \xi) \leq C m(\d \xi)$, the following estimates hold: there exists a constant $C$ depending only on $\alpha,d$, s.t. $\forall v\in \R_+\backslash \{0\}, z\in \R^d$,
\begin{equation}\label{hold-ARONSON_STABLE}
	C^{-1}v^{-\frac{d}{\alpha}}\left( 1+ \frac{|z|}{v^{\frac{1}{\alpha}}} \right)^{-(d+\alpha)}\leq p_\alpha(v,z)\leq Cv^{-\frac{d}{\alpha}}\left( 1+ \frac{|z|}{v^{\frac{1}{\alpha}}} \right)^{-(d+\alpha)}.
\end{equation}

On the other hand let us mention that the sole non-degeneracy condition  
\begin{equation*}
	\kappa^{-1} |\lambda|^\alpha \leq \int_{\mathbb{S}^{d-1}} |\lambda \cdot \xi|^\alpha \mu(\d \xi) \leq \kappa |\lambda|^\alpha,
\end{equation*} 
does not allow to derive \textit{global} heat kernel estimates for the noise density. In \cite{Wa07}, Watanabe investigates the behavior of the density of an $\alpha$-stable process in terms of properties fulfilled by the support of its spectral measure $\mu$. From this work, we know that whenever the measure $\mu$ is not equivalent to the Lebesgue measure $m$ on the unit sphere, accurate estimates on the density of the stable process are delicate to obtain. \\


From now on, and in particular in Section \ref{hold-SEC_PROOF_REG} which is dedicated to the proof of technical lemmas, we will be using the \textit{proxy} notation 
\begin{equation}
	\bar{p}_\alpha (v,z) := \begin{cases}
		C_\alpha {v^{-\frac{d}{\alpha}}} \left(1+\frac{|z|}{v^{\frac{1}{\alpha}}} \right)^{-(d+\alpha)} &\mbox{ if } \alpha \in (1,2)\\
		{(2\pi c v)^{-\frac{d}{2}}}  \exp \left( -c^{-1}\frac{|z|^2}{2v}\right), \ c \ge 1 &\mbox{ if } \alpha = 2,
	\end{cases},\qquad\qquad v>0,z\in \R^d,\label{hold-DEF_P_BAR}
\end{equation}
where, for $\alpha\in (1,2)$, $C_\alpha$ is chosen so that $\forall v>0, \int \bar{p}_\alpha (v,y) \d y = 1$, and $c:=c(d)$
is a global given constant for $\alpha=2$. We will explicitly rely on the global bounds provided by $\bar{p}_\alpha$. 
Observe importantly that,	 keeping in mind that \eqref{hold-ARONSON_STABLE}, in the pure jump case, there exists $C\ge 1$ s.t. for all $(v,z)\in \R_+\backslash\{0\}\times \R^d $,
\begin{equation}\label{hold-EQ_PALPHA_BARPALPHA_JUMPS}
	C^{-1}\bar p_\alpha(v,z)\le  p_\alpha(v,z) \le C\bar p_\alpha(v,z),
\end{equation}
and the results could be stated with either the proxy density $\bar p_\alpha$ or the density $p_\alpha $ of the noise itself. However, the equivalence in \eqref{hold-EQ_PALPHA_BARPALPHA_JUMPS} fails in the Gaussian case, due to the exponential tails. This is why the results will be stated in terms of $\bar p_\alpha $. Observe as well that from the definition in \eqref{hold-DEF_P_BAR} we readily have the following important properties:
\begin{trivlist}
	\item[-] (Approximate) convolution property: there exists a constant ${\mathfrak c}\ge 1$ s.t. for all $u,v\in \R_+\backslash\{0\},\   x,y\in \R^d$,
	\begin{equation}
		\label{hold-APPR_CONV_PROP}
		\int_{\R^d} \bar p_\alpha(u,z-x)\bar p_\alpha(v,y-z) \d z\le {\mathfrak c} \bar p_\alpha(u+v,y-x). 
	\end{equation}
	In particular, for $\alpha=2 $ the convolution is \textit{exact} and ${\mathfrak c}=1$.
	\item[-] Time-scale for the spatial moments: for all $0\le \delta<\alpha, \alpha\in (1,2) $ and for all $\delta\ge 0 $ if $\alpha=2 $, there exists $C_{\alpha,\delta} $ s.t.
	\begin{equation}
		\label{hold-spatial-moments}
		\int_{\R^d} |z|^\delta\bar p_\alpha(v,z)\d z\le C_{\alpha,\delta} v^{\frac{\delta}{\alpha}}.
	\end{equation}
	
\end{trivlist}

Further properties related to the density of the driving noise, notably concerning its {time-space derivatives}, are stated in Lemma \ref{hold-lemma-stable-sensitivities} below.
\subsection{Main results}\label{hold-subsec-main-result}
We first give some important estimates concerning the densities of the SDE \eqref{hold-sde} and its associated Euler scheme \eqref{hold-scheme-interpo}.

\begin{propo}[Density estimates for the diffusion and its Euler scheme]\label{hold-prop-main-estimates}
	The unique weak solution to Equation \eqref{hold-sde} starting from $x$ at time $s\in [0,T]$ admits for all $t\in (s,T]$  a density $\Gamma(s,x,t,\cdot) $. Furthermore there exists a constant $C:=C({d},b,\alpha,T)$ s.t.
	for all $y\in \R^d $ the following upper-bound holds:
	\begin{align}\label{hold-ineq-density-diff}
		&\Gamma(s,x,t,y)\le C \bar p_\alpha (t-s,y-x),
	\end{align}
	with $\bar p_\alpha $ defined in \eqref{hold-DEF_P_BAR}, as well as the following control for the H\"older regularity in the forward time variable:
	\begin{equation}\label{hold-holder-time-gamma}
		\forall 0\le s<t<t'\le T,\    |t-t'|\le (t-s) ,\   |\Gamma(s,x,t,y)-\Gamma (s,x,t',y)| \le C \frac{(t'-t)^\frac{\gamma}{\alpha}}{(t-s)^\frac{\gamma}{\alpha}} \bar p_\alpha (t'-s,y-x).
	\end{equation}
	Also, {for $\varepsilon \in (0,{\gamma\wedge 1}] $, there exists $C_\varepsilon :=C_\varepsilon({d},b,\alpha,T)$} s.t. forall $0\le s<t\le T,\ x,y,w\in \R^d$ s.t.  $|y-w|\le (t-s)^{\frac 1\alpha} $,
	\begin{align}
		\label{hold-holder-space-gamma}
		|\Gamma(s,x,t,y)-\Gamma(s,x,t,w)|\le {C_\varepsilon}\left(\frac{|y-w|}{(t-s)^{\frac1\alpha}} \right)^{{\gamma_\varepsilon}}\bar p_\alpha(t-s,w-x),\ {\gamma_\varepsilon:=(\gamma\wedge 1)-\varepsilon}.
	\end{align}

	Similarly, for any positive integer $n$ and $h=\frac{T}{n}$, the corresponding Euler scheme $X^h$  defined  in \eqref{hold-scheme-interpo} starting from $x\in \R^d$ at time $t_k:=kh, k\in \llbracket 0,n-1 \rrbracket$ admits for $t\in (t_k,T]$ a  transition 
	density $\Gamma^h(t_k,x,t,\cdot) $, for which, for all $y\in \R^d $:
	\begin{align}\label{hold-ineq-density-scheme}
		&\Gamma^h(t_k,x,t,y)\le C \bar p_\alpha (t-t_k,y-x).
	\end{align}
	Also, for all $0< t_j<t_k\le T $, $x,y,w\in \R^d $, $|y-w|\le (t_k-t_j)^{\frac 1\alpha} $,
	\begin{align}
		\label{hold-holder-space-gammah}
		|\Gamma^h(t_j,x,t_k,y)-\Gamma^h(t_j,x,t_k,w)|\le {C_\varepsilon} \left(\frac{|y-w|{+h^{\frac 1\alpha}}}{(t_k-t_j)^{\frac1 \alpha}} \right)^{{\gamma_\varepsilon}}\bar p_\alpha(t_k-t_j,w-x).
	\end{align}
\end{propo}   
Existence of the densities and the related Aronson type bounds \eqref{hold-ineq-density-diff} and \eqref{hold-ineq-density-scheme} readily follow from
\cite{JM211} and \cite{FJM24}. The sensitivity controls \eqref{hold-holder-time-gamma}, \eqref{hold-holder-space-gamma} and \eqref{hold-holder-space-gammah} are proven in Section \ref{hold-SEC_PROOF_REG}.

\begin{remark}[About additional controls on the density of the SDE and the Euler scheme]
	Let us point out that the densities $\Gamma,\Gamma^h $ also satisfy additional controls. Namely, some gradient controls in the backward spatial variable could be established. Anyhow, for the error analysis, these controls are not needed. They could be derived following the approach of \cite{FJM24}.
\end{remark}	

The main result of the chapter is then the following theorem.    
\begin{theorem}[Convergence Rate for the stable-driven Euler scheme with $L_t^\infty \mathcal{C}_x^{\beta}$ drift]\label{hold-thm-main}
	Denoting by $\Gamma $ and $\Gamma^h$ the respective densities of the SDE \eqref{hold-sde} and its Euler scheme defined in \eqref{hold-euler-scheme}, there exists a constant $C:=C({d},b,\alpha,T)<\infty$ s.t. for  all $h=T/n$ with $n\in\N^*$, and all $t\in(0,T]$, $x,y\in \R^d $,
	\begin{align}
		|\Gamma^h(0,x,t,y)-\Gamma(0,x,t,y)| &\leq C\big(1+t^{-\frac \beta\alpha}\big)  h^{\frac{\gamma}{\alpha}} \bar p_\alpha(t,y-x),\label{hold-BD_THM}
	\end{align}
	where $\gamma=\beta+\alpha-1 >0$ is again the ``gap to singularity" defined in \eqref{hold-serrin}.
\end{theorem}
\begin{remark}[Weak error involving an additional test function]
	Let us mention that if one is interested in the weak error for some test function $f$,
	$\mathcal E(f,x,t,h):=\E_{0,x}[f(X_t^h)-f(X_t)],$
	as soon as $f$ is $\delta\in [\beta,1] $-H\"older (not necessarily bounded) then, a rate can be derived as a consequence of the convergence of $|\Gamma(s,x,t,y)-\Gamma^h(s,x,t,y)|$ using a simple cancellation argument:
	\begin{align*}
		\mathcal E(f,x,t,h)=&\int_{\R^d}(\Gamma^h-\Gamma)(0,x,t,y)f(y)\d y=\int_{\R^d}(\Gamma^h-\Gamma)(0,x,t,y)\big(f(y)-f(x)\big)\d y,\\
		|\mathcal E(f,x,t,h)|\le& C h^{\frac \gamma\alpha}\big(1+t^{-\frac \beta\alpha}\big)\int_{\R^d} \bar p_\alpha(t,y-x) |x-y|^\delta\d y\underset{\eqref{hold-DEF_P_BAR}}{\le} \tilde C \big(1+t^{-\frac \beta\alpha}\big) t^{\frac\delta \alpha}h^{\frac \gamma\alpha}.
	\end{align*}	
	{Precisely, the smoothness of $f$ allows to absorb the time-singularity from \eqref{hold-BD_THM} in small time.}
\end{remark}

\section{Proof of the main results}   \label{hold-SEC_PROOF_THM}

We begin this section with recalling some quantitative properties of the density of the driving noise as well as the Duhamel representations of the densities which will be the starting point to analyze the corresponding error.

\subsection{Representation and Estimates on the densities of the diffusion and its  Euler scheme}\label{hold-subsec-scheme-estimates}
As in the papers \cite{JM211}, \cite{FJM24} in which the weak error was investigated for Lebesgue drifts, we will expand the densities of the SDE and its Euler scheme along the underlying heat equation. In particular, since the drift we consider is here bounded, existence of the density and related Aronson type upper-bounds readily follow from these works (see e.g.  Propositions 2.1. and 2.3 in \cite{JM211} for the Brownian case and Theorem 1 and Proposition 1 in \cite{FJM24} for the pure jump one).\\

In the current H\"older setting in space, in order to explicitly take advantage of the additional spatial regularity we will rely on appropriate cancellation techniques. We start recalling some useful controls for the underlying heat kernel $p_\alpha $, the density of the driving noise $Z$.

\subsubsection{Controls on the density of the stable noise} \label{hold-subsec-controls-palpha}
\begin{propo}[Density estimates for the heat equation]\label{hold-prop-controls-palpha} There exists $C\ge 1,\bar c\ge 1$ s.t. for all $0 \le s<t\le T $, $(x,w) \in (\R^d)^2$, and any muti-index $\zeta \in \N^d $, $|\zeta|\in \{1, 2\} $, 
	$\theta\in \{0,1\} $,
	\begin{align}
		|\partial_t^\theta\nabla_x^\zeta p_\alpha (t-s,w-x)|\le C (t-s)^{-(\theta+\frac{|\zeta|+d}{\alpha})}\left(1+\frac{\left| w-x\right|}{(t-s)^{\frac{1}{\alpha}}} \right)^{-d-\alpha-|\zeta|},\ \alpha\in (1,2),\notag\\
		|\partial_t^\theta\nabla_x^\zeta p_\alpha (t-s,w-x)|\le C(t-s)^{-(\theta+\frac{|\zeta|+d}{\alpha})}g_{\bar c}(t-s,w-x),\ \alpha=2, \label{hold-GD_BOUNDS}
	\end{align}
	where $g_{\bar c}(u,z):= \frac{1}{(2\pi \bar c)^{\frac d2}}\exp(-\frac{|z|^2}{2\bar c u})$ stands for the Gaussian density of the Gaussian vector with variance $\bar cI_d $.
	In particular, under \eqref{hold-serrin}, 
	\begin{align}\label{hold-drift-smoothing-iso-noise}
		|w-x|^\beta|\partial_t^\theta\nabla_{{x}}^\zeta  p_\alpha(t-s,w-x)|\le C (t-s)^{-\theta+\frac{\beta-|\zeta|}{\alpha}} \bar p_\alpha(t-s,w-x),
	\end{align}
	and, for $w'\in \R^d$ s.t. $|w-w'|\lesssim (t-s)^{\frac{1}{\alpha}}$, 
	\begin{align}\label{hold-drift-smoothing-iso-noise-diag}
		|w-x|^\beta|\partial_t^\theta\nabla_x^\zeta p_\alpha\left(t-s,w-x+(w-w')\right)|\le C (t-s)^{-\theta+\frac{\beta-|\zeta|}{\alpha}} \bar p_\alpha(t-s,w-x).
	\end{align}
\end{propo}
\begin{proof}
	The estimates in \eqref{hold-GD_BOUNDS} are plain to prove directly if $\alpha=2$. Turning now to $\alpha\in (1,2) $, since we have assumed $Z$ to be isotropic, it is well known that the {spatial} derivatives of the density of the driving noise enjoy better concentration properties, see e.g. Lemma 2.8  in \cite{MZ22}, which proves \eqref{hold-GD_BOUNDS} for $\theta=0 $. For $\theta=1$, let us proceed as follows: for $|w-x|\leq (t-s)^{\frac{1}{\alpha}}$, the bound follows from the Fourier representation of the density:
	\begin{equation}
		p_\alpha (t,x)=(2\pi)^{-d}\int \exp \left[i x \cdot \xi\right] \exp \left[-c_\alpha t |\xi|^\alpha \right] \d \xi, \qquad c_\alpha>0.
	\end{equation}
	For  $|w-x|\geq (t-s)^{\frac{1}{\alpha}}$, let us recall that the density of the isotropic stable process can be expressed as:
	\begin{align*}
		p_\alpha(t-s,w-x)&=\int_0^\infty g(r,w-x) p_{S^{\frac \alpha2}}\left(t-s, r\right)\d r\\
		&=\int_0^\infty g(r,w-x) \frac 1{(t-s)^{\frac 2\alpha}}p_{S^{\frac \alpha2}}\left(1,\frac r{(t-s)^{\frac2\alpha}}\right)\d r,
	\end{align*}
	where $g:(r,x)\in\R_+\times \R^d \mapsto (2\pi r)^{-d/2} \exp \left[-|x|^2/(2r)\right]$ denotes the standard gaussian density and $p_{S^{\alpha/2}}$ stands for the density of the $\alpha/2$ stable subordinator. Hence, integrating by parts, 
	\begin{align*}
		\partial_t p_\alpha(t-s,w-x)&=\int_0^\infty g(r,w-x) \partial_t\left(\frac 1{(t-s)^{\frac 2\alpha}}p_{S^{\alpha/2}}\left(1,\frac r{(t-s)^{\frac2\alpha}}\right)\right) \d r\\
		&=-\frac{2}{\alpha}\Bigg[\frac{1}{t-s}\int_0^\infty g(r,w-x) \frac 1{(t-s)^{\frac 2\alpha}}p_{S^{\frac \alpha2}}\left(1,\frac r{(t-s)^{\frac2\alpha}}\right)\d r\\ & \qquad \qquad +\int_0^{\infty} g(r,w-x)\frac{1}{(t-s)^{\frac 2\alpha}} \partial_r p_{S^{\alpha/2}}\left(1,\frac r{(t-s)^{\frac2\alpha}}\right) \frac{r}{t-s} \d r\Bigg]\\
		&=-\frac{2}{\alpha(t-s)}\Bigg[p(t-s,w-x)-\int_0^{\infty} \partial_r (r g(r,w-x)) \frac{1}{(t-s)^{\frac 2\alpha}} p_{S^{\alpha/2}}\left(1,\frac r{(t-s)^{\frac2\alpha}}\right) \d r\Bigg] \\
		&=-\frac{2}{\alpha (t-s)} \int_0^{\infty}r \partial_r g(r,w-x) \frac{1}{(t-s)^{\frac 2\alpha}} p_{S^{\alpha/2}}\left(1,\frac r{(t-s)^{\frac2\alpha}}\right) \d r.
	\end{align*}	
	Recalling that $	\partial_r g(r,x)= \Big(-\frac{d}{2r} +\frac{|x|^2}{2 r^2}\Big) g(r,x)$, we have, for $ \zeta$ s.t. $ |\zeta|\in \{1,2\} $,
	\begin{equation}
		|r\nabla_x^\zeta\partial_r g(r,x)| \lesssim \left(\frac{1}{r^{\frac{|\zeta|}{2}}} + \left(\frac{|x|}{r}\right)^{|\zeta|}\right) r^{-\frac{d}{2}}\exp \left[-\frac{|x|^2}{2r}\right]\lesssim\frac{1}{r^{\frac{|\zeta|+d}{2}}} \exp \left[-\lambda\frac{|x|^2}{2r}\right],
	\end{equation}
	for some $\lambda>1$. Plugging this into the previous equation and setting $u=|w-x|^2/r$, we get
	\begin{align*}
		&|\nabla_x^\zeta \partial_t p_\alpha(t-s,w-x)|\lesssim \frac{1}{t-s}\int_0^{\infty} \frac{1}{r^{\frac{|\zeta|+d}{2}}} \exp \left[-\lambda\frac{|w-x|^2}{2r}\right] \frac{1}{(t-s)^{\frac 2\alpha}} p_{S^{\alpha/2}}\left(1,\frac r{(t-s)^{\frac2\alpha}}\right) \d r\\
		&\qquad \lesssim \frac{1}{t-s}\int_0^{\infty} \left(\frac{|w-x|^2}{u}\right)^{-\frac{|\zeta|+d}{2}} \exp \left[-\frac{\lambda u}{2}\right] \frac{1}{(t-s)^{\frac 2\alpha}} p_{S^{\alpha/2}}\left(1,\frac{|w-x|^2}{u(t-s)^{\frac2\alpha}}\right) \frac{|w-x|^2}{u^2}\d u.
	\end{align*}
	Using now the global bound on the law of the stable subordinator (see e.g. the proof of Lemma 4.8 in \cite{KKM16}, in particular how the global bound is derived from Equation (4.8)):
	$$\forall s>0, \qquad p_{S^{\alpha/2}}(1,s)\lesssim s^{-1-\frac{\alpha}{2}}, $$
	we have,
	\begin{align*}
		|\nabla_x^\zeta \partial_t p_\alpha(t-s,w-x)|&\lesssim \frac{1}{t-s}\int_0^{\infty} \left(\frac{|w-x|^2}{u}\right)^{-\frac{|\zeta|+d}{2}} \exp \left[-\frac{\lambda u}{2}\right] \frac{1}{(t-s)^{\frac 2\alpha}} \left(\frac{|w-x|^2}{u(t-s)^{\frac2\alpha}}\right)^{-1-\frac{\alpha}{2}} \frac{|w-x|^2}{u^2}\d u\\
		&\lesssim |w-x|^{-|\zeta|-d-\alpha}\int_0^\infty u^{\frac{|\zeta|+d+\alpha}{2}-1}\exp \left[-\frac{\lambda u}{2}\right] \d u\lesssim |w-x|^{-|\zeta|-d-\alpha}\\
		&\lesssim ((t-s)^{\frac{1}{\alpha}}+|w-x|)^{-|\zeta|-d-\alpha} \lesssim (t-s)^{-\big(1+\frac{|\zeta|+d}{\alpha}\big)}\left(1+\frac{\left| w-x \right|}{(t-s)^{\frac{1}{\alpha}}} \right)^{-d-\alpha-|\zeta|},
	\end{align*}
	recalling for the last line that $|w-x|>(t-s)^{\frac{1}{\alpha}}$. This concludes the proof of \eqref{hold-GD_BOUNDS}.\\

		An important consequence of the above estimates is precisely \eqref{hold-drift-smoothing-iso-noise}. Namely, for  $\alpha\in (1,2) $ we get 
		\begin{align*}
			|w-x|^\beta|\partial_t^\theta\nabla_x^\zeta  p_\alpha (t-s,w-x)| &\le C (t-s)^{-\theta+\frac{\beta-|\zeta|}{\alpha}}\left(\frac{|w-x|}{(t-s)^\frac{1}{\alpha}}\right)^\beta \left(1+\frac{\left| w-x\right|}{(t-s)^{\frac{1}{\alpha}}} \right)^{-|\zeta|}  \bar  p_\alpha (t-s,w-x)\\
			&\le C (t-s)^{-\theta+\frac{\beta-|\zeta|}{\alpha}} \bar  p_\alpha (t-s,w-x).
		\end{align*}
		Also, for $\alpha=2 $, 
		\begin{align*}
			|w-x|^\beta|\partial_t^\theta\nabla_x^\zeta p_\alpha (t-s,w-x)| &\le C (t-s)^{-\theta+\frac{\beta-|\zeta|}2}\left( \frac{|w-x|}{(t-s)^{\frac 12}}\right)^\beta g_{\bar c}(t-s,w-x)\\
			&\le C (t-s)^{-\theta+\frac{\beta-|\zeta|}{2}} \bar  p_2 (t-s,w-x).
		\end{align*}   
		In that case the concentration constant is slightly deteriorated  whereas in the pure jump case we took advantage of the concentration improvement for the derivatives. In any case \eqref{hold-drift-smoothing-iso-noise} is proven. {Equation \eqref{hold-drift-smoothing-iso-noise-diag} follows from the same proof as \eqref{hold-drift-smoothing-iso-noise} noting that in the diagonal regime $|w-w'|\leq (t-s)^{\frac{1}{\alpha}}$,
			\begin{align*}
				\bar p_\alpha( t-s,w-x+(w-w'))&\le C \bar p_\alpha( t-s,w-x),\ \alpha\in (1,2),\\
				g_{\bar c}(t-s,w-x+(w-w'))& \le C g_{\tilde  c}(t-s,w-x), \tilde c<c,\ \alpha=2.
			\end{align*}
			Hence, $|w-w'|$ can be seen as a negligible perturbation.}
	\end{proof}

	\subsubsection{Duhamel representation for the densities}   
	
	To compute the error rate, we will start from the following Duhamel representations, which are proved respectively in \cite{JM211} and \cite{FJM24} for $\alpha=2 $ and $\alpha\in (1,2) $:
	\begin{propo}[Duhamel representations for the densities of the SDE and the Euler scheme]\label{hold-prop-main-estimates-D}
		The density $\Gamma(s,x,t,\cdot) $ of the unique weak solution to Equation \eqref{hold-sde} starting from $x$ at time $s\in [0,T)$ admits {the} following  Duhamel representation: for all ${t\in (s, T],\  y\in \R^d }$,
		\begin{align}
			\Gamma(s,x,t,y)
			= p_\alpha(t-s,y-x)-\int_{s}^{ t}\E_{s,x}\left[b(r,X_r)\cdot\nabla_y  p_\alpha(t-r,y-X_r)\right]\d r,\label{hold-duhamel-Diff}
		\end{align} 
		where the expectation subscript means that $X_s=x$.\\
		Similarly, for $k\in \llbracket 0,n-1\rrbracket ,\ t\in (t_k,T] $, the density of $X_t^h$ admits, conditionally to $X_{t_k}^h=x$, a transition density $\Gamma^h(t_k,x,t,\cdot) $, which again enjoys a Duhamel type representation: for all $y\in \R^d $,
		\begin{align}
			\Gamma^h(t_k,x,t,y)
			= p_\alpha(t-t_k,y-x)-\int_{t_k}^{ t}\E_{t_k,x}\left[b(U_{\tau_r^h/h},X^h_{\tau_r^h})\cdot\nabla_y  p_\alpha(t-r,y-X^h_r)\right]\d r,\label{hold-duhamel-scheme}
		\end{align} 
		where the expectation subscript means that $X_{t_k}^h=x$.
	\end{propo}
	\subsection{Proof of Theorem \ref{hold-thm-main}}
	In this section we will use for two quantities $A$ and $B$ the symbol $A\lesssim B $ whenever there exists a constant $C:=C({d},b,\alpha,T)$ s.t. $A\le CB $. Namely,
	\begin{equation}\label{hold-DEF_LESSSIM}
		A\lesssim B \Longleftrightarrow \exists C:=C({d},b,\alpha,T),\ A\le CB.
	\end{equation}

	Starting from Proposition \ref{hold-prop-main-estimates-D} and  comparing the Duhamel formula of the scheme, \eqref{hold-duhamel-scheme}, to that of the diffusion, \eqref{hold-duhamel-Diff}, we get
	\begin{align*}
		&	\Gamma^h(0,x,t,y)-\Gamma(0,x,t,y)\\ & \qquad =\E_{0,x} \left[ \int_0^t \left( b(s,X_s)\cdot \nabla_y p_\alpha(t-s,y-X_s)-b(U_{\tau_s^h/h},X_{\tau_s^h}^h)\cdot\nabla_y p_\alpha(t-s,y-X_s^h)\right) \d s \right].
	\end{align*}
	In the previous equation and the rest of this chapter, we denote $\E_{0,x}[\cdot]:=\E[\cdot|X_0=X_0^h=x]$.
	We will split the error in the following way: 
	\begin{align}
		&\Gamma^h(0,x,t,y)-\Gamma(0,x,t,y)\notag\\
		& \qquad =\int_0^h \E_{0,x} \left[ b(s,X_s)\cdot \nabla_y p_\alpha (t-s,y-X_s)-b(U_0,{x})\cdot\nabla_y p_\alpha (t-s,y-X_s^h) \right] \d s\notag\\
		& \qquad \qquad + \int_h^{\tau_t^h-h} \E_{0,x} \left[ b(s,X_s)\cdot \nabla_y p_\alpha(t-s,y-X_s)-b(s,X_{\tau_s^h})\cdot \nabla_y  p_\alpha (t-s,y-X_{\tau_s^h})
		\right] \d s\notag\\
		& \qquad \qquad + \int_h^{\tau_t^h-h} \E_{0,x} \bigg[ b(s,X_{\tau_s^h})\cdot \nabla_y p_\alpha (t-s,y-X_{\tau_s^h})
		-b(s,X_{\tau_s^h}^h)\cdot \nabla_y  p_\alpha (t-s,y-X_{\tau_s^h}^h) \bigg] \d s\notag\\
		&\qquad \qquad + \int_h^{\tau_t^h-h} \E_{0,x} \left[b(U_{\tau_s^h/h},X_{\tau_s^h}^h)\cdot \left( \nabla_y p_\alpha(t-s,y-X_{\tau_s^h}^h)-\nabla_y  p_\alpha (t-s,y-X_s^h)\right)\right]\d s\notag\\
		& \qquad \qquad+\int_h^{\tau_t^h-h} \E_{0,x} \left[b(U_{\tau_s^h/h},X_{\tau_s^h}^h)\cdot \left( \nabla_y p_\alpha(t-U_{\tau_s^h/h},y-X_{\tau_s^h}^h)-\nabla_y  p_\alpha (t-s,y-X_{\tau_s^h}^h)\right)\right]\d s\notag\\
		&\qquad\qquad+\int_{\tau_t^h-h}^t \E_{0,x} \left[ b(s,X_s)\cdot \nabla_y p_\alpha (t-s,y-X_s)-b(U_{\tau_s^h/h},X_{\tau_s^h}^h)\cdot\nabla_y p_\alpha (t-s,y-X_s^h)  \right] \d s\notag\\
		&\qquad =:\Delta_1 + \Delta_2 + \Delta_3 + \Delta_4 + \Delta_5+\Delta_6,\label{hold-DECOUP_ERR}
	\end{align}
	where we exploited that:
	\begin{align}
		&\int_h^{\tau_t^h-h} \E_{0,x} \left[b(U_{[\tau_{s}^h/h]},X_{\tau_s^h}^h)\cdot \nabla_y  p_\alpha (t-U_{[\tau_s^h/h]},y-X_{\tau_s^h}^h)\right]\d s\nonumber\\
		&\qquad \qquad =\sum_{i=1}^{\tau_t^h/h-1}\frac{1}{h}\int_{t_i}^{t_{i+1}} \int_{t_i}^{t_{i+1}}
		\E_{0,x} \left[b(r,X_{t_i}^h)\cdot \nabla_y  p_\alpha (t-r,y-X_{t_i}^h)\right] \d s \d r\nonumber\\
		&\qquad \qquad =\int_h^{\tau_t^h{-h}} \E_{0,x} \left[b(s,X_{\tau_s^h}^h)\cdot \nabla_y  p_\alpha (t-s,y-X_{\tau_s^h}^h)\right]\d s,\label{hold-convenient-fubini}
	\end{align}
	for the correspondence between the last term in $\Delta_3$ and the first one in $\Delta_5 $.\\

	For $\Delta_1$, we rely on the fact that we work on the first time step, and thus we do not even need the smoothing effect in $(t-s)^\beta$ provided by the drift in \eqref{hold-drift-smoothing-iso-noise}. Let us first expand the expectation:
	\begin{align*}
		\Delta_1 &=  \int_0^h \int \Big(\Gamma(0,x,s,z)b(s,z)\cdot \nabla_y  p_\alpha (t-s,y-z) - \Gamma^h(0,x,s,z)b(U_0,x) \cdot \nabla_y p_\alpha (t-s,y-z)\Big)\d z\d s.
	\end{align*}
	Assuming w.l.o.g. that $t> 2h$ so that for $s\in (0,h)$, $(t-s)^{-\frac{1}{\alpha}}<(t/2)^{-\frac{1}{\alpha}}$, then using \eqref{hold-ineq-density-scheme}, \eqref{hold-ineq-density-diff}, 
	\eqref{hold-GD_BOUNDS} and the boundedness of $b$, we have
	\begin{align}
		|\Delta_1|&\lesssim \int_0^h \int \frac{1}{(t-s)^{\frac{1}{\alpha}}}\bar p_\alpha (s,z-x) \bar p_\alpha (t-s,y-z)  \d z\d s\nonumber\\
		&\lesssim \bar p_\alpha (t,y-x)\int_0^h \frac{1}{(t-s)^{\frac{1}{\alpha}}} \d s\lesssim \bar  p_\alpha (t,y-x)h t^{-\frac 1\alpha}\nonumber\\
		&\lesssim \bar p_\alpha (t,y-x)h^{\frac{\gamma}{\alpha}}h^{\frac{1-\beta}{\alpha}}t^{-\frac 1\alpha}\lesssim \bar p_\alpha (t,y-x)h^{\frac{\gamma}{\alpha}}t^{-\frac \beta\alpha},\label{hold-maj-D1-mainthm}
	\end{align}
	{exploiting the \textit{convolution property} \eqref{hold-APPR_CONV_PROP} of the density $\bar p_\alpha $\footnote{which is actually just an \textit{approximate} convolution property in the pure jump case.} for the second inequality
		and} recalling that $h\le t$ for the last inequality. \\
	
	Let us turn to $\Delta_2$. Expanding the inner expectation and using the time regularity of $\Gamma$ in the forward time variable (see \eqref{hold-holder-time-gamma})
	\begin{align*}
		|\Delta_2|&=\left| \int_h^{\tau_t^h-h} \int\left[\Gamma (0,x,s,z)-\Gamma (0,x,\tau_s^h,z)\right]  b(s,z)\cdot\nabla_y  p_\alpha (t-s,y-z) \d z \d s\right|\\
		&\lesssim \int_h^{t}  \int\frac{(s-\tau_s^h)^{\frac{\gamma}{\alpha}}}{(\tau_s^h)^{\frac{\gamma}{\alpha}}} \bar p_\alpha(s,z-x) (t-s)^{-\frac{1}{\alpha}}\bar p_\alpha (t-s,y-z)\d z \d s,
	\end{align*}
	using as well {\eqref{hold-GD_BOUNDS}} for the last inequality. {Again, from} \eqref{hold-APPR_CONV_PROP},  along with the fact that $s-\tau_s^h\leq h$ and that for $s\geq h, (\tau_s^h)^{-1}\leq 2s^{-1}$, we can write
	\begin{align}
		|\Delta_2| &\lesssim \bar p_\alpha (t,y-x) h^{\frac{\gamma}{\alpha}}  \int_h^{t}  s^{-\frac{\gamma}{\alpha}}(t-s)^{{-\frac{1}{\alpha}}} \d s\lesssim \bar p_\alpha (t,y-x) h^{\frac{\gamma}{\alpha}}. \label{hold-maj-D2-mainthm}
	\end{align}
	
	The term $\Delta_3$, which is the one that will allow to apply a Gronwall type argument, will be treated at the end of the current error analysis.\\
	
	Let us turn to the term $\Delta_4$ in \eqref{hold-DECOUP_ERR}, which is by far the more delicate.
	Let us then introduce the following lemma:
	\begin{lemma}[Smoothing effect of the drift]\label{hold-THE_LEMMA_FOR_THM} Let $\zeta$ be a multi-index with length $1\leq |\zeta|\leq {2}$ and $\delta\in \{0,1\}$. Then, for all $(x,\mathfrak{y})\in (\R^d)^2$, $0\leq h\leq t/2 \leq s < \tau_t^h-h \leq T$, $r>0$,
		\begin{align}
			&\left|\int \Gamma^{h}(0,x,\tau_s^h,z)b(r,z) \partial_t^\delta\nabla_{\mathfrak{y}}^\zeta p_\alpha (t-{\tau_s^h},\mathfrak{y}-z) \d z\right|	\lesssim \bar{p}_\alpha(t,\mathfrak{y}-x)(t-\tau_s^h)^{-\delta + \frac{\beta-|\zeta|}{\alpha}}\left(1+(\tau_s^h)^{-\frac{\beta}{\alpha}}\right)\label{hold-lemma-time-gain-thm}.
		\end{align}
	\end{lemma} 
	\begin{proof} Let use the following cancellation argument:
		\begin{align*}
			&I:= \int \Gamma^{h}(0,x,\tau_s^h,z)b(r,z) \partial_t^\delta \nabla_{\mathfrak{y}}^\zeta p_\alpha (t-{\tau_s^h},\mathfrak{y}-z) \d z\\
			&\qquad =\int [\Gamma^{h}(0,x,\tau_s^h,z)b(r,z)-\Gamma^{h}(0,x,\tau_s^h,\mathfrak{y})b(r,\mathfrak{y})] \partial_t^\delta\nabla_{\mathfrak{y}}^\zeta p_\alpha (t-{\tau_s^h},\mathfrak{y}-z) \d z.
		\end{align*}
		Then, {using the regularity of $b$ and \eqref{hold-holder-space-gammah}, taking therein $\varepsilon $ such that $\gamma_\varepsilon\ge \beta $}, we have
		\begin{align}
			&|\Gamma^{h}(0,x,\tau_s^h,z)b(r,z)-\Gamma^{h}(0,x,\tau_s^h,\mathfrak{y})b(r,\mathfrak{y})|\notag\\
			&\qquad \lesssim \bar{p}_\alpha (\tau_s^h,z-x)|\mathfrak{y}-z|^\beta + \Vert b \Vert_{L^\infty} \frac{(|\mathfrak{y}-z| + h^{\frac{1}{\alpha}})^{\beta}}{(\tau_s^h)^{\frac{\beta}{\alpha}}}\left[\bar{p}_\alpha(\tau_s^h,z-x)+\bar{p}_\alpha(\tau_s^h,\mathfrak{y}-x)\mathbb{1}_{|\mathfrak{y}-z|^\alpha\geq \tau_s^h}\right]\label{hold-USEFUL_FOR_DELTA_6}\\
			&\qquad \lesssim \bar{p}_\alpha (\tau_s^h,z-x)|\mathfrak{y}-z|^\beta  \left(1+(\tau_s^h)^{-\frac{\beta}{\alpha}}\right)+  \frac{(|\mathfrak{y}-z|^\beta+ (t-\tau_s^h)^{\frac{\beta}{\alpha}})}{(\tau_s^h)^{\frac{\beta}{\alpha}}}\bar{p}_\alpha(t,\mathfrak{y}-x),\notag
		\end{align}
		where we used the fact that we consider times $s\geq t/2$ to write $\bar{p}_\alpha(\tau_s^h,\mathfrak{y}-x)\lesssim\bar{p}_\alpha(t,\mathfrak{y}-x)$ and the fact that for $s\leq \tau_t^h-h$, $h\leq t-\tau_s^h$ {for the last inequality}. Plugging this into $I$ and using \eqref{hold-drift-smoothing-iso-noise}, we get \eqref{hold-lemma-time-gain-thm}.
	\end{proof}
	Going back to the bound for $\Delta_4$, conditioning w.r.t. $\sigma(X_{\tau_s}^h,U_{\tau_s^h/h}) $ (the sigma-algebra generated by the random variables $X_{\tau_s^h}$ and $U_{\tau_s^h}$) and using the harmonicity of the (gradient of the) stable heat kernel (or Itô's formula between $\tau_s^h $ and $s$) in order to get rid of the noise increment, we can write that for any bounded and measurable $\varphi :\R^d\times \R \rightarrow \R^d$,
	\begin{equation}
		\E\left[\nabla_y p_\alpha \left(t-s,y-\varphi(X_{\tau_s^h}^h,U_{\tau_s^h})+(Z_s-Z_{\tau_s^h})\right)|\sigma(X_{\tau_s}^h,U_{\tau_s^h/h}) \right]=\nabla_y p_\alpha \left(t-\tau_s^h,y-\varphi(X_{\tau_s^h}^h,U_{\tau_s^h})\right).
	\end{equation}
	Using this, we get
	\begin{align}
		\Delta_4&=\int_h^{t/2} \E_{0,x} \left[b(U_{\tau_s^h/h},X_{\tau_s^h}^h)\cdot \left( \nabla_y p_\alpha(t-s,y-X_{\tau_s^h}^h)-\nabla_y  p_\alpha (t-s,y-X_s^h)\right)\right]\d s\notag\\
		&\qquad +\int_{t/2}^{\tau_t^h-h} \E_{0,x} \left[b(U_{\tau_s^h/h},X_{\tau_s^h}^h)\cdot \left( \nabla_y p_\alpha(t-s,y-X_{\tau_s^h}^h) -\nabla_y  p_\alpha (t-s,y-(X_{\tau_s^h}^h+Z_s-Z_{\tau_s^h}))\right)\right]\d s\notag\\
		&\qquad +\int_{t/2}^{\tau_t^h-h} \E_{0,x} \left[b(U_{\tau_s^h/h},X_{\tau_s^h}^h)\cdot\left( \nabla_y  p_\alpha (t-s,y-(X_{\tau_s^h}^h+Z_s-Z_{\tau_s^h}))\right.\right.\notag\\
		&\quad\qquad \qquad\qquad\qquad \qquad\left.\left.-\nabla_y  p_\alpha (t-s,y-(X_{\tau_s^h}^h+b(U_{\tau_s^h/h},X_{\tau_s^h}^h)(s-\tau_s^h)+Z_s-Z_{\tau_s^h}))\right)\right] \d s\nonumber\\
		&=\int_h^{t/2} \E_{0,x} \left[b(U_{\tau_s^h/h},X_{\tau_s^h}^h)\cdot \left( \nabla_y p_\alpha(t-s,y-X_{\tau_s^h}^h)-\nabla_y  p_\alpha (t-{\tau_s^h},y-{(X_{\tau_s^h}^h+b(U_{\tau_s^h/h},X_{\tau_s^h}^h)(s-\tau_s^h))})\right)\right]\d s\notag\\
		&\qquad +\int_{t/2}^{\tau_t^h-h} \E_{0,x} \left[b(U_{\tau_s^h/h},X_{\tau_s^h}^h)\cdot \left( \nabla_y p_\alpha(t-s,y-X_{\tau_s^h}^h)- \nabla_y  p_\alpha (t-\tau_s^h,y-X_{\tau_s^h}^h)\right)\right]\d s\notag\\
		&\qquad +\int_{t/2}^{\tau_t^h-h} \E_{0,x} \left[b(U_{\tau_s^h/h},X_{\tau_s^h}^h)\cdot\left( \nabla_y  p_\alpha (t-\tau_s^h,y-X_{\tau_s^h}^h)\right.\right.\notag\\
		&\quad\qquad \qquad\qquad\qquad \qquad\left.\left.-\nabla_y  p_\alpha (t-\tau_s^h,y-(X_{\tau_s^h}^h+b(U_{\tau_s^h/h},X_{\tau_s^h}^h)(s-\tau_s^h)))\right)\right] \d s\nonumber\\ &=:\Delta_{41}+\Delta_{42}+\Delta_{43}.\label{hold-NEW_DECOUP_DELTA_4}
	\end{align}   
	For $\Delta_{41}$, there is no need to compensate for singularities in $(t-s)$ on the considered time interval:
	\begin{align}
		|\Delta_{41}|& =\nonumber \left|\int_h^{t/2} \frac{1}{h}\int_{\tau_s^h}^{\tau_s^h+h}\int \Gamma^h(0,x,\tau_s^h,z)b(r,z)\big[\nabla_y p_\alpha(t-s,y-z)- \nabla_y  p_\alpha (t-\tau_s^h,y-z)\right. \\ & \qquad \qquad \left. + \nabla_y p_\alpha(t-\tau_s^h,y-z)- \nabla_y  p_\alpha (t-\tau_s^h,y-z-b(r,z)(s-\tau_s^h))\big]\d z \d r \d s\right|\nonumber\\
		& \lesssim \int_h^{t/2} \int \bar{p}_\alpha(\tau_s^h,z-x)\Vert b \Vert_{L^\infty} \left[\frac{(s-\tau_s^h)}{(t-s)^{1+\frac{1}{\alpha}}} +\frac{(s-\tau_s^h)\Vert b \Vert_{L^\infty}}{(t-\tau_s^h)^{\frac{2}{\alpha}}} \right]\bar p_\alpha(t-\tau_s^h,y-z)\d z \d s\nonumber\\
		& \lesssim \bar p_\alpha(t,y-x)h^{\frac{\gamma}{\alpha}}\int_h^{t/2} \left[(t-s)^{-\frac{\gamma+1}{\alpha}}+(t-\tau_s^h)^{1-\frac{\gamma+2}{\alpha}}\right]\d s \nonumber\\
		&\lesssim \bar p_\alpha(t,y-x)h^{\frac{\gamma}{\alpha}}t^{1-\frac{\gamma+1}{\alpha}}= \bar p_\alpha(t,y-x)h^{\frac{\gamma}{\alpha}}t^{-\frac{\beta}{\alpha}}.\label{hold-CTR_DELTA_41}
	\end{align}
	For $\Delta_{42}$, we will make the same time sensitivity appear and then expand it with a Taylor formula:
	\begin{align}
		\Delta_{42}& = \int_{t/2}^{\tau_t^h-h} \frac{1}{h}\int_{\tau_s^h}^{\tau_s^h+h}\int \Gamma^h(0,x,\tau_s^h,z)b(r,z)\cdot \left[\nabla_y p_\alpha(t-s,y-z)- \nabla_y  p_\alpha (t-\tau_s^h,y-z)\right]\d z \d r \d s\nonumber\\
		&= {-}\int_{t/2}^{\tau_t^h-h} \frac{1}{h}\int_{\tau_s^h}^{\tau_s^h+h}\int \Gamma^h(0,x,\tau_s^h,z)b(r,z)\cdot \int_0^1 \partial_t \nabla_yp_\alpha(t-s+\mu(s-\tau_s^h),y-z)(s-\tau_s^h)\d \mu\d z \d r \d s{.}\nonumber
	\end{align}
	Then, using the same {cancellation} techniques as in the proof of \eqref{hold-lemma-time-gain-thm} with $\mathfrak{y}=y$, we get
	\begin{align}
		|\Delta_{42}| &\lesssim {(1+t^{-\frac{\beta}{\alpha}})}\int_{t/2}^{\tau_t^h-h}  \int_0^1 \frac{(s-\tau_s^h)}{(t-s+\mu(s-\tau_s^h))^{1+{\frac{1-\beta}{\alpha}}}} \bar p_\alpha(t+\tau_s^h-s+\mu(s-\tau_s^h),y-{x})\d \mu \d s\nonumber\\
		&\lesssim \bar p_\alpha(t,y-{x})(1+t^{-\frac{\beta}{\alpha}})\int_{t/2}^{\tau_t^h-h} \frac{(s-\tau_s^h)}{(t-\tau_s^h)^{1+\frac{1-\beta}{\alpha}}}  \d s\nonumber\\
		&\lesssim \bar p_\alpha(t,y-{x})(1+t^{-\frac{\beta}{\alpha}}){h}\int_{t/2}^{\tau_t^h-h}{(t-s)^{-1+\frac{\beta-1}{\alpha}}  }\d s\nonumber\\
		&\lesssim \bar p_\alpha(t,y-{x})(1+t^{-\frac{\beta}{\alpha}}){h(t-\tau_t^h+h)^{\frac{\beta-1}{\alpha}}}\lesssim \bar p_\alpha(t,y-z)(1+t^{-\frac{\beta}{\alpha}})h^{\frac{\gamma}{\alpha}}\label{hold-CTR_DELTA_42}.
	\end{align}
	Let us turn to $\Delta_{43}$, let us write
	\begin{align}
		\Delta_{43}&=\int_{t/2}^{\tau_t^h-h}\frac{1}{h}\int_{\tau_s^h}^{\tau_s^h+h}\int \Gamma^h(0,x,\tau_s^h,z) b(r,z)\cdot\big[ \nabla_y  p_\alpha (t-\tau_s^h,y-z)-{\nabla_w  p_\alpha (t-\tau_s^h,w)|_{w=y-z-b(r,y)(s-\tau_s^h)}}\notag\\
		&\quad\qquad \qquad+{\nabla_w  p_\alpha (t-\tau_s^h,w)|_{w=y-z-b(r,y)(s-\tau_s^h)}}-\nabla_y  p_\alpha (t-\tau_s^h,y-(z+b(r,z)(s-\tau_s^h)))\big]\d z \d r \d s\nonumber\\
		&:=\Delta_{431}+\Delta_{432}.\label{hold-DECOUP_DELTA_43}
	\end{align}
	{We carefully mention that this additional pivot is needed in order to use cancellation arguments for the first term (in order to have a drift which does not depend on the spatial integration variable) and to take full force of the regularity of the drift for the second one}.
	For $\Delta_{431}$, we use a Taylor expansion and then \eqref{hold-lemma-time-gain-thm}:
	\begin{align}
		|\Delta_{431}|&=\bigg|\int_{t/2}^{\tau_t^h-h}\frac{1}{h}\int_{\tau_s^h}^{\tau_s^h+h}\int \Gamma^h(0,x,\tau_s^h,z) \nonumber b(r,z)\\ & \qquad \qquad \qquad\cdot \int_0^1 {\nabla_w^2  p_\alpha (t-\tau_s^h,w)|_{w=y-z-\mu b(r,y)(s-\tau_s^h)}}\cdot  b(r,y)(s-\tau_s^h)\d \mu \d z \d r \d s\bigg|\nonumber\\
		&\lesssim (1+t^{-\frac{\beta}{\alpha}}) \int_{t/2}^{\tau_t^h-h} {\frac 1h\int_{\tau_s^h}^{\tau_s^h+h}}\int_0^1\bar  p_\alpha (t,y-\mu b(r,y)(s-\tau_s^h)-x)  \Vert b \Vert_{L^\infty} \frac{(s-\tau_s^h)}{(t-\tau_s^h)^{\frac{2-\beta}{\alpha}}}\d \mu {\d r} \d s\nonumber\\	
		&\lesssim (1+t^{-\frac{\beta}{\alpha}})\bar  p_\alpha (t,y-x)h  \int_{t/2}^{\tau_t^h-h} (t-\tau_s^h)^{-\frac{2-\beta}{\alpha}}\d s.\nonumber
	\end{align}
	{Recalling that $\bar  p_\alpha (t,y-\mu b(r,y)(s-\tau_s^h)-x)\lesssim \bar  p_\alpha (t,y-x) $ for the second inequality (with a slight notational abuse in the gaussian case since the variance is then modified).
		Now,	 if $(2-\beta)/\alpha <1 \iff \alpha+\beta>2 $ (which is e.g. always the case in the Brownian setting), the time integral is convergent and uniformly bounded in $h$. The term $\Delta_{431} $ then has order $h$ regarding the time step. If now $(2-\beta)/\alpha =1 $, it has order $h|\ln(h)|\lesssim h^{\frac \gamma\alpha} $.  We can thus assume w.l.o.g. that $(2-\beta)/\alpha >1 $. Then,}
	\begin{align}	
		|\Delta_{431}|	&\lesssim (1+t^{-\frac{\beta}{\alpha}}) \bar  p_\alpha (t,y-x)h (t-\tau_t^h+h)^{1-\frac{2-\beta}{\alpha}}\nonumber\\
		&\lesssim  (1+t^{-\frac{\beta}{\alpha}}) \bar  p_\alpha (t,y-x)h^{\frac{\gamma}{\alpha}+(1-\frac{1}{\alpha})}\lesssim \bar  p_\alpha (t,y-x) (1+t^{-\frac{\beta}{\alpha}}) h^{\frac{\gamma}{\alpha}},\label{hold-CTR_DELTA_431}
	\end{align}
	
	Let us turn to $\Delta_{432}$, which we first expand with a Taylor formula:
	\begin{align*}
		&\Delta_{432}=\int_{t/2}^{\tau_t^h-h}\frac{1}{h}\int_{\tau_s^h}^{\tau_s^h+h}\int \Gamma^h(0,x,\tau_s^h,z) \nonumber b(r,z)\\ &\qquad  \cdot \int_0^1{ \nabla_w^2  p_\alpha \left(t-\tau_s^h,w\right)|_{w=y-z-\left[b(r,z)-\mu \left[b(r,y)-b(r,z)\right]\right](s-\tau_s^h)}}\cdot  [b(r,y)-b(r,z)](s-\tau_s^h)\d \mu \d z \d r \d s.\nonumber
	\end{align*}
	Since $t-\tau_s^h \geq s-\tau_s^h$, using the boundedness and the H\"older regularity of $b$ and then \eqref{hold-drift-smoothing-iso-noise-diag}, we have
	\begin{align}
		&{\left|\nabla_w^2  p_\alpha \left(t-\tau_s^h,w\right)|_{w=y-z-\left[b(r,z)-\mu \left[b(r,y)-b(r,z)\right]\right](s-\tau_s^h)}\cdot  [b(r,y)-b(r,z)]\right|}\nonumber\\
		&\qquad \lesssim (t-\tau_s^h)^{\frac{\beta-2}{\alpha}}\bar p_\alpha(t-\tau_s^h, y-z),\label{hold-KILLING_THE_DRIFT-thm}
	\end{align}
	thus yielding, with the same computations of the time integrals as for $\Delta_{431}$,
	\begin{align}
		{|}\Delta_{432}{|}&\lesssim {(1+t^{-\frac \beta\alpha})}\bar{p}_\alpha(t,y-x)\int_{t/2}^{\tau_t^h-h} (s-\tau_s^h)(t-\tau_s^h)^{\frac{\beta-2}{\alpha}}\d s\nonumber\\
		&\lesssim {(1+t^{-\frac \beta\alpha})}\bar{p}_\alpha(t,y-x)h^{\frac{\gamma}{\alpha}}.\label{hold-CTR_DELTA_432}
	\end{align}
	Plugging \eqref{hold-CTR_DELTA_431} and \eqref{hold-CTR_DELTA_432} into \eqref{hold-DECOUP_DELTA_43} and 
	from \eqref{hold-CTR_DELTA_41}, \eqref{hold-CTR_DELTA_42} and \eqref{hold-NEW_DECOUP_DELTA_4} , we obtain
	\begin{equation}\label{hold-maj-D4-mainthm}
		|\Delta_4|\lesssim (1+t^{-\frac \beta \alpha})\bar p_\alpha(t,y-x) h^{\frac \gamma\alpha}.\\ 
	\end{equation}
	
	Observe now that the term $\Delta_5$ in \eqref{hold-DECOUP_ERR} could be handled just as $\Delta_{41}$ { and $\Delta_{42}$} 
	above. Indeed, the time variable is there randomized, but once expanded through the density, the quantity is really similar. {Importantly, it makes a pure time sensitivity of the stable kernel appear}. This therefore yields
	\begin{equation}\label{hold-maj-D5-mainthm}
		|\Delta_5|\lesssim (1+t^{-\frac \beta \alpha})\bar p_\alpha(t,y-x) h^{\frac \gamma\alpha}. \\
	\end{equation}    
	
	It thus remains to handle the contribution $\Delta_6 $ associated with the last time steps. The quantities involved can actually again be estimated {using cancellation arguments}. 
%
	Write:
	\begin{align}
		|\Delta_6|&\le  \left|\int_{\tau_t^h-h}^t \E_{0,x} \left[ b(s,X_s)\cdot \nabla_y p_\alpha (t-s,y-X_s)-b(U_{\tau_s^h/h},X_{\tau_s^h}^h)\cdot\nabla_y p_\alpha (t-s,y-X_s^h)  \right] \d s\right|\notag\\
		&\le  \left|\int_{\tau_t^h-h}^t  \int_{}\Big(\Gamma(0,s,x,z)b(s,z)-\Gamma(0,s,x,y)b(s,y)\Big) \nabla_y p_\alpha (t-s,y-z) \d z \d s\right|\notag\\
		&\quad +  
		\left|\int_{\tau_t^h-h}^t \E_{0,x} \left[ b(U_{\tau_s^h/h},X_{\tau_s^h}^h)\cdot\nabla_y p_\alpha (t-s,y-X_s^h)  \right] \d s\right|\notag\\
		&=:\Delta_{61}+\Delta_{62}.\label{DECOUP_DELTA_6}
	\end{align}
	Using the Hölder estimate \eqref{hold-holder-space-gamma} for the density of the diffusion (still taking therein $\varepsilon$ s.t. $\gamma_\varepsilon\ge \beta $), and exploiting \eqref{hold-drift-smoothing-iso-noise} as in the proof of Lemma \eqref{hold-lemma-time-gain-thm}, we then get:
	\begin{equation}\label{hold-maj-D61-mainthm}
		|\Delta_{61}|\lesssim \int_{\tau_t^h-h}^t  \Big((t-s)^{-\frac 1\alpha+\frac\beta\alpha}+{(t-s)^{-\frac 1\alpha}h^{\frac \beta\alpha}}\Big)(1+(\tau_s^h)^{-\frac{\beta}\alpha})\bar p_\alpha (t,y-x)\d s\lesssim(1+t^{-\frac \beta \alpha})\bar p_\alpha(t,y-x) h^{\frac \gamma\alpha}. 
	\end{equation}       
Because of the time difference appearing for the contribution of the scheme, i.e. $\Delta_{62} $ involves both $X_{\tau_s^h}^h $ and $X_s^h $, the corresponding term requires a more careful analysis. 

Namely, the previous cancellation argument cannot be reproduced exactly.

Write:
 \begin{align}
 \Delta_{62}
:=&\int_{\tau_t^h-h}^t \E_{0,x} \left[ b(U_{\tau_s^h/h},X_{\tau_s^h}^h)\cdot\nabla_y p_\alpha (t-s,y-X_s^h)  \right] \d s\notag\\
=&\int_{\tau_t^h-h}^t \E_{0,x} \left[\Big( b(U_{\tau_s^h/h},X_{\tau_s^h}^h)\cdot\nabla_y p_\alpha \big(t-\tau_s^h,y-(X_{\tau_s^h}^h+b(U_{\tau_s^h/h},X_{\tau_s^h}^h)(s-\tau_s^h)) \big) \right] \d s\notag\\
=&\int_{\tau_t^h-h}^t \E_{0,x} \left[\Big( b(U_{\tau_s^h/h},X_{\tau_s^h}^h)\cdot\nabla_y p_\alpha (t-\tau_s^h,y-X_{\tau_s^h}^h)  \right] \d s\notag\\
&+\int_{\tau_t^h-h}^t \E_{0,x} \left[\Big( b(U_{\tau_s^h/h},X_{\tau_s^h}^h)\cdot \right.\notag\\
&\hspace*{2cm}\left.\Big(\nabla_y p_\alpha \big(t-\tau_s^h,y-(X_{\tau_s^h}^h+b(U_{\tau_s^h/h},X_{\tau_s^h}^h)(s-\tau_s^h))\big)-\nabla_y p_\alpha \big(t-\tau_s^h,y-X_{\tau_s^h}^h\big)\Big)  \right] \d s.\notag\\
=:& \Delta_{621}+ \Delta_{622},\label{NEW_SPLIT_LT_OLD}
 \end{align}
 where we again used the harmonicity of the stable heat kernel for the first equality. 
 The term $\Delta_{621} $ can be handled as the previous term $\Delta_{61} $. Namely,
 \begin{align*}
&|\Delta_{621}|\\
\le& \int_{\tau_t^h-h}^t \frac{\d r}{h} \left|\int_{\tau_t^h-h}^t  \int_{}\Big(\Gamma^h(0,\tau_s^h,x,z)b(r,z)-\Gamma^h(0,\tau_s^h,x,y)b(r,y)\Big) \nabla_y p_\alpha (t-\tau_s^h,y-z) \d z \d s\right|\\
\le& C\int_{\tau_t^h-h}^t  \int_{} \Big(\bar p_\alpha(\tau_s^h,y-z)+\bar p_\alpha(\tau_s^h,y-x)\Big)(1+(\tau_{s}^h)^{-\frac \beta\alpha})\Big(\frac{1}{(t-s)^{\frac 1\alpha-\frac \beta\alpha}}+\frac{h^{\beta}}{(t-s)^{\frac1\alpha}}\Big) \bar p_\alpha(t-\tau_s^h,y-z)\d z\\
\le& Ch^{\frac{\alpha-1+\beta}\alpha}\bar p_\alpha(t,y-x),
 \end{align*}
where we used \eqref{hold-holder-space-gammah} (taking therein $\gamma_\varepsilon=\beta $) for the second inequality. Also, for $\theta\in (0,\alpha-1) $,
 \begin{align*}
|\Delta_{622}|&\le \int_{\tau_t^h-h}^t \bar p_\alpha(\tau_s^h,z-x)\frac{(\|b\|_\infty h)^\theta}{(t-\tau_s^h)^{\frac 1\alpha+\frac{\theta}{\alpha}}}\bar p_\alpha(t-\tau_s^h,y-z) \d s\\
&\le C\bar p_\alpha(t,y-x) h^{\frac{\alpha-(1+\theta)}{\alpha}+\theta} =C\bar p_\alpha(t,y-x) h^{\frac{\alpha-1}{\alpha}+\frac{\theta(\alpha-1)}{\alpha}}.
\end{align*}
Hence, this approach works provided one can take $\theta(\alpha-1)=\beta\iff \theta=\frac{\beta}{\alpha-1} $. The integrability constraint $\theta<\alpha-1 $ then yields $\beta<(\alpha-1)^2 $. Hence, this will always work in the Brownian setting, but not necessarily in the strictly stable one.
	
Write now in whole generality:		
 \begin{align}
 \Delta_{62}
:=&\int_{\tau_t^h-h}^t \E_{0,x} \left[ b(U_{\tau_s^h/h},X_{\tau_s^h}^h)\cdot\nabla_y p_\alpha (t-s,y-X_s^h)  \right] \d s\notag\\
=&\int_{\tau_t^h-h}^t \E_{0,x} \left[\Big( b(U_{\tau_s^h/h},X_{\tau_s^h}^h)-b(U_{\tau_s^h/h},X_{s}^h)\Big)\cdot\nabla_y p_\alpha (t-s,y-X_s^h)  \right] \d s\notag\\
&+\int_{\tau_t^h-h}^t \E_{0,x} \left[ b(U_{\tau_s^h/h},X_{s}^h)\cdot\nabla_y p_\alpha (t-s,y-X_s^h)  \right] \d s.\notag\\
=:& \Delta_{621}+ \Delta_{622}.\label{NEW_SPLIT_LT}
 \end{align}
Then:
 \begin{align}
&\Delta_{621}\notag\\
=&\int_{\tau_t^h-h}^t \E_{0,x} \left[\Big( b(U_{\tau_s^h/h},X_{\tau_s^h}^h)-b(U_{\tau_s^h/h},X_{s}^h)\Big)\cdot\nabla_y p_\alpha (t-s,y-(X_{\tau_s^h}^h+b(X_{\tau_s^h}^h)(s-\tau_s^h)+Z_{s}-Z_{\tau_s^h}))  \right] \d s\notag
\\
=&\int_{\tau_t^h-h}^t\frac{1}{h}\int_{\tau_s^h}^{\tau_s^h+h}\E_{0,x} \Big[\Big( b(r,X_{\tau_s^h}^h)-b(r,X_{\tau_s^h}^h+b(r,X_{\tau_s^h}^h)(s-\tau_s^h)+Z_s-Z_{\tau_s^h}))\Big)\notag\\
&\cdot\nabla_y p_\alpha (t-s,y-(X_{\tau_s^h}^h+b(X_{\tau_s^h}^h)(s-\tau_s^h)+Z_{s}-Z_{\tau_s^h}))  \Big] \d r\d s\notag\\
=&\int_{\tau_t^h-h}^t\frac{1}{h}\int_{\tau_s^h}^{\tau_s^h+h}\E_{0,x} \Big[ \E_{\mathcal F_{\tau_s^h}}[\mathcal E(s,r,X_{\tau_s^h}^h,Z_s-Z_{\tau_s^h})]\Big]\d r \d s.
\label{PREAL_DELTA_621}
 \end{align}
 We now consider the conditional expectation:
 \begin{align}
&\E_{\mathcal F_{\tau_s^h}}[\mathcal E(s,r,X_{\tau_s^h}^h,Z_s-Z_{\tau_s^h})]\notag\\
=&\int_{\R^d} \Big( b(r,X_{\tau_s^h}^h)-b(r,X_{\tau_s^h}^h+b(r,X_{\tau_s^h}^h)(s-\tau_s^h)+z))\Big)\nabla_y p_\alpha (t-s,y-(X_{\tau_s^h}^h+b(X_{\tau_s^h}^h)(s-\tau_s^h)+z))\notag\\
&\times p_\alpha(s-\tau_s^h ,z) \d z\notag\\
=&-\int_{\R^d} \Big( b(r,X_{\tau_s^h}^h)-b(r,X_{\tau_s^h}^h+b(r,X_{\tau_s^h}^h)(s-\tau_s^h)+z))\Big)\nabla_z p_\alpha (t-s,y-(X_{\tau_s^h}^h+b(X_{\tau_s^h}^h)(s-\tau_s^h)+z))\notag\\
&\times p_\alpha(s-\tau_s^h ,z) \d z\notag\\
=&\int_{\R^d} \Big( b(r,X_{\tau_s^h}^h)-b(r,X_{\tau_s^h}^h+b(r,X_{\tau_s^h}^h)(s-\tau_s^h)+z))\Big) p_\alpha (t-s,y-(X_{\tau_s^h}^h+b(X_{\tau_s^h}^h)(s-\tau_s^h)+z))\notag\\
&\times \nabla_zp_\alpha(s-\tau_s^h ,z) \d z\notag\\
&+\int_{\R^d} {\rm div}_z\big(b(r,X_{\tau_s^h}^h+b(r,X_{\tau_s^h}^h)(s-\tau_s^h)+z)\big)p_\alpha (t-s,y-(X_{\tau_s^h}^h+b(X_{\tau_s^h}^h)(s-\tau_s^h)+z))p_\alpha(s-\tau_s^h ,z) \d z\notag\\
=:&\Delta_{6211}+\Delta_{6212}. \label{EC_LAST_STEPS}
 \end{align}
  One then gets:
  \begin{align*}
|\Delta_{6211}|\lesssim \int_{\R^d}(|z|^\beta+h^{\beta})p_\alpha (t-s,y-(X_{\tau_s^h}^h+b(X_{\tau_s^h}^h)(s-\tau_s^h)+z)) |\nabla_zp_\alpha(s-\tau_s^h ,z)| \d z.
  \end{align*}
 Observing that for any stability index $\alpha\in (1,2] $,  from \eqref{hold-GD_BOUNDS}:
 $$|z|^\beta |\nabla_zp_\alpha(s-\tau_s^h ,z)|\lesssim (s-\tau_s^h)^{-\frac{1-\beta}\alpha}\bar p_\alpha(s-\tau_s^h,z),$$
 which is indeed obvious in the Gaussian case and follows from the improved concentration of the tails for the derivatives of the isotropic stable heat kernel, we actually get:
   \begin{align}
|\Delta_{6211}|&\lesssim \big( h^\beta(s-\tau_s^h)^{-\frac 1\alpha}+(s-\tau_s^h)^{-\frac{1-\beta}\alpha}\big)\bar p_\alpha(t-\tau_s^h,y-(X_{\tau_s^h}^h+b(X_{\tau_s^h}^h)(s-\tau_s^h)))\notag\\
&\lesssim \big( h^\beta(s-\tau_s^h)^{-\frac 1\alpha}+(s-\tau_s^h)^{-\frac{1-\beta}\alpha}\big)\bar p_\alpha(t-\tau_s^h,y-X_{\tau_s^h}^h),\label{CTR_I6211}
  \end{align}
noting again that, since $s-\tau_s^h\le t-\tau_s^h $, the remaining local transition in the first inequality can indeed be removed.
 On the other hand, noting that ${\rm div}_z\big(b(r,\cdot)\big)\in B_{\infty,\infty}^{\beta-1}$ we get by duality, see \eqref{dual-ineq}, and observing that $\|{\rm div}_z\big(b(r,\cdot)\big)\|_{B_{\infty,\infty}^{\beta-1}} \lesssim \|b(r,\cdot)\|_{B_{\infty,\infty}^\beta}$:
   \begin{align}
|\Delta_{6212}|&\lesssim \|b(r,\cdot)\|_{B_{\infty,\infty}^\beta}\|p_\alpha (t-s,y-(X_{\tau_s^h}^h+b(X_{\tau_s^h}^h)(s-\tau_s^h)+\cdot))p_\alpha(s-\tau_s^h ,\cdot)\|_{B_{1,1}^{1-\beta}}\notag\\
&\lesssim \|b(r,\cdot)\|_{B_{\infty,\infty}^\beta}\|p_\alpha (t-s,y-(X_{\tau_s^h}^h+b(X_{\tau_s^h}^h)(s-\tau_s^h)+\cdot))p_\alpha(s-\tau_s^h ,\cdot)\|_{B_{1,1}^{1-\beta}}\notag\\
&\lesssim \|b(r,\cdot)\|_{B_{\infty,\infty}^\beta} p_\alpha (t-\tau_s^h,y-(X_{\tau_s^h}^h+b(X_{\tau_s^h}^h)(s-\tau_s^h)))(t-\tau_s^h)^{-\frac{1-\beta}{\alpha}}\Big(1+\frac{(t-\tau_s^h)^{\frac \zeta\alpha}}{(t-s)^{\frac \zeta\alpha}}+\frac{(t-\tau_s^h)^{\frac \zeta\alpha}}{(s-\tau_s^h)^{\frac \zeta\alpha}}\Big)\notag\\
&\lesssim \|b(r,\cdot)\|_{B_{\infty,\infty}^\beta} p_\alpha (t-\tau_s^h,y-X_{\tau_s^h}^h)(t-\tau_s^h)^{-\frac{1-\beta}{\alpha}}\Big(1+\frac{(t-\tau_s^h)^{\frac \zeta\alpha}}{(t-s)^{\frac \zeta\alpha}}+\frac{(t-\tau_s^h)^{\frac \zeta\alpha}}{(s-\tau_s^h)^{\frac \zeta\alpha}}\Big)\label{CTR_I6212}
 \end{align}
for $\zeta\in(1-\beta,1) $, using  \eqref{BESOV_DENSITY_LP} for the last bound. Plugging \eqref{CTR_I6211} and \eqref{CTR_I6212} into \eqref{EC_LAST_STEPS} and \eqref{PREAL_DELTA_621}, we eventually derive:
\begin{align}
&|\Delta_{621}\notag|\notag\\
\lesssim& \int_{\tau_t^h-h}^{t}\Big[ h^\beta(s-\tau_s^h)^{-\frac 1\alpha}+(s-\tau_s^h)^{-\frac{1-\beta}\alpha}+(t-\tau_s^h)^{-\frac{1-\beta}{\alpha}}\Big(1+\frac{(t-\tau_s^h)^{\frac \zeta\alpha}}{(t-s)^{\frac \zeta\alpha}}+\frac{(t-\tau_s^h)^{\frac \zeta\alpha}}{(s-\tau_s^h)^{\frac \zeta\alpha}}\Big)\Big]\notag\\
&\times 
\int_{\R^d} \bar p_\alpha(0,x,\tau_s^h,w)\bar p_\alpha(t-\tau_s^h,y-w) \d w \d s\notag\\
\lesssim& \bar p_\alpha(t,y-x) h^{\frac{\alpha-1+\beta}{\alpha}},\label{CTR_DELTA_621}
\end{align}
which is the expected bound.

%

 The term $ \Delta_{622}  $ in \eqref{NEW_SPLIT_LT} can be handled as the previous $\Delta_{61} $ (same cancellation argument, using the estimates established in the proof of \eqref{hold-holder-space-gammah}, where the estimate is shown for an arbitrary final time and not only a discretization time) and would as well lead to the indicated control as in \eqref{hold-maj-D61-mainthm}. Hence, from \eqref{CTR_DELTA_621}, \eqref{hold-maj-D61-mainthm} and \eqref{DECOUP_DELTA_6}, we eventually derive:
\begin{equation}\label{hold-maj-D6-mainthm}
|\Delta_6|\lesssim \int_{\tau_t^h-h}^t  \Big((t-s)^{-\frac 1\alpha+\frac\beta\alpha}+\textcolor{black}{(t-s)^{-\frac 1\alpha}h^{\frac \beta\alpha}}\Big)(1+(\tau_s^h)^{-\frac{\beta}\alpha})\bar p_\alpha (t,y-x)\d s\lesssim(1+t^{-\frac \beta \alpha})\bar p_\alpha(t,y-x) h^{\frac \gamma\alpha}. 
\end{equation}       

	Gathering estimates \eqref{hold-maj-D1-mainthm}, \eqref{hold-maj-D2-mainthm}, \eqref{hold-maj-D4-mainthm}, \eqref{hold-maj-D5-mainthm}, \eqref{hold-maj-D6-mainthm}, we have 
	\begin{align}
		&|\Gamma^h(0,x,t,y)-\Gamma (0,x,t,y)|\lesssim \bar p_\alpha(t,y-x) h^{\frac{\gamma}{\alpha}}(1+t^{-\frac \beta\alpha})\nonumber \\
		& \qquad + \left| \int_h^{t} \E_{0,x} \bigg[ b(s,X_{\tau_s^h})\cdot \nabla_y p_\alpha (t-\tau_s^h,y-X_{\tau_s^h})-b(s,X_{\tau_s^h}^h)\cdot \nabla_y p_\alpha (t-\tau_s^h,y-X_{\tau_s^h}^h) \bigg] \d s \right|.
	\end{align}
	Set 
	\begin{equation}
		f_h(u) = \sup_{(x,z)\in (\R^d)^2} \frac{|\Gamma^h(0,x,u,z)-\Gamma (0,x,u,z)|}{\bar p_\alpha(u,z-x)}.
	\end{equation}
	Observe from \eqref{hold-ineq-density-scheme} and \eqref{hold-ineq-density-diff} that $f_h$ is bounded uniformly in $h$ and the time variable. We then have, using \eqref{hold-drift-smoothing-iso-noise}, the boundedness of $b$ and the convolution property of the stable kernel,
	\begin{align*}
		\frac{|\Gamma^h(0,x,t,y)-\Gamma (0,x,t,y)|}{p_\alpha (t,y-x)}&\lesssim h^{\frac{\gamma}{\alpha}}(1+t^{-\frac \beta\alpha})+  \frac{1}{\bar p_\alpha (t,y-x)}\int_h^{t}\frac{f_h(\tau_s^h)}{(t-\tau_s^h)^{\frac{1}{\alpha}}}\int  \bar p_\alpha(\tau_s^h,z-x)\bar p_\alpha (t-\tau_s^h,y-z)\d z\d s\\
		&\lesssim h^{\frac{\gamma}{\alpha}}(1+t^{-\frac \beta\alpha})+ \int_h^{t}\frac{f_h(\tau_s^h)}{(t-\tau_s^h)^{\frac{1}{\alpha}}}\d s.
	\end{align*}
	The previous bound being uniform in $x$ and $y$, we get
	\begin{align*}
		f_h(t) &\lesssim h^{\frac{\gamma}{\alpha}}(1+t^{-\frac \beta\alpha})+ \int_h^{t}\frac{f_h(\tau_s^h)}{(t-\tau_s^h)^{\frac{1}{\alpha}}}\d s.
	\end{align*}
	Using a \textit{discrete} Gr\"onwall-Volterra lemma, we obtain Theorem \ref{hold-thm-main}.

	\section{Proof of the regularity results from Proposition \ref{hold-prop-main-estimates}}
	
	\label{hold-SEC_PROOF_REG}
	{ This section is devoted to the proof of the controls \eqref{hold-holder-space-gamma}, \eqref{hold-holder-space-gammah} concerning the Hölder continuity in space in the forward variable for the density of the diffusion and the scheme respectively as well as to \eqref{hold-holder-time-gamma}, Hölder regularity in time in the forward time variable. These estimates were crucial in order to prove the main theorem in the previous section. Importantly, we achieve respectively the order $\gamma_\varepsilon=(\gamma\wedge 1)-\varepsilon,\ \varepsilon>0 $ in space and $\gamma/\alpha $ in time.}\\
	
	
	{   Note that the lower exponent attained in space  is sufficient for   the previous proof of the main result to work.  We actually used the spatial regularity of the densities of the diffusion and the scheme in a cancellation argument involving a product with the drift (see former \eqref{hold-USEFUL_FOR_DELTA_6}). In this context, we just need in practice the lower order $\beta $-regularity corresponding to the spatial smoothness of the drift. Indeed, for $\beta\in(0,1) $ one can find $\varepsilon $ s.t. $(\gamma\wedge 1)-\varepsilon\ge \beta $. We mention as well that this first spatial estimate actually allows in a second time to derive the expected exponent $\gamma $ when $\gamma<1 $ for the H\"older regularity in the forward variable. The proof is provided for completeness (since we actually insisted on the parabolic bootstrap phenomenon) for the diffusion in Appendix \ref{hold-BOOT}.  
	}\\

	We start this section recalling some usual yet important controls on the density of the driving noise that we will profusely use in order to prove \eqref{hold-holder-time-gamma}. The proof is somehow standard and can be e.g. found in \cite{FJM24} for $\alpha\in (1,2) $ and \cite{JM211} in the gaussian case.
	\begin{lemma}[Stable sensitivities - Estimates on the $\alpha $-stable kernel]\label{hold-lemma-stable-sensitivities}
		For each multi-index  $\zeta$ with length $|\zeta|\leq {2}$, and for all $0<u\leq u'\leq T$, $(x,x')\in (\R^d)^2$, $\theta\in (0,1] $,
		\begin{itemize} 
			\item Time H\"older regularity: 
			\begin{align}\label{hold-holder-time-palpha}
				\left|\nabla_x^\zeta p_\alpha (u,x)-\nabla_x^\zeta p_\alpha (u',x)\right| \lesssim \frac{|u-u'|^\theta}{u^{\theta+\frac{|\zeta|}{\alpha}}} \left(\bar p_\alpha (u,x)+\bar p_\alpha (u',x)\right).
			\end{align}
			\item Spatial H\"older regularity: 
			\begin{align}\label{hold-holder-space-palpha}
				\left|\nabla_x^\zeta p_\alpha (u,x)-\nabla_{x'}^\zeta p_\alpha (u,x')\right| \lesssim \left(\frac{|x-x'|^\theta}{u^{\frac{\theta}{\alpha}}} \wedge 1\right)\frac{1}{u^{\frac{|\zeta|}{\alpha}}}\left(\bar p_\alpha (u,x)+\bar p_\alpha (u,x')\right).
			\end{align}
		\end{itemize}
	\end{lemma}

	\subsection{Proof of the spatial regularity}
	
	We start this section providing the estimate \eqref{hold-holder-space-gamma} for the diffusion. We could actually have established \eqref{hold-holder-space-gammah} only and then derive \eqref{hold-holder-space-gamma} passing to the limit exploiting the convergence in law of the Euler scheme to the diffusion, which would have allowed to transfer the estimates on densities. However, we provide a complete proof on the diffusion first since it is actually simpler than the one for the scheme and already emphasizes the key ideas, which will as well appear in the proof of \eqref{hold-holder-time-gamma} (time regularity).
	
	\subsubsection{Proof of \eqref{hold-holder-space-gamma}: forward spatial H\"older regularity for the diffusion}
	
	Define, for $\eta >0 $ meant to be small (viewed as a spatial viscosity parameter),
	$$h_{s,x}^{\varepsilon,\eta}(t):=\sup_{(z,z')\in (\R^d)^2} \left\{ \frac{|\Gamma(s,x,t,z)-\Gamma(s,x,t,z')|(t-s)^{\frac{\gamma_\varepsilon}{\alpha}}}{\left(\bar p_\alpha(t-s,z-x)+\bar p_\alpha(t-s,z'-x)\right){(|z-z'| \vee \eta)^{\gamma_\eps}}} \right\}.$$
	{Since we already know from \eqref{hold-ineq-density-diff} that $ \frac{|\Gamma(s,x,t,z)-\Gamma(s,x,t,z')|}{\bar p_\alpha(t-s,z-x)+\bar p_\alpha(t-s,z'-x)}<\infty$}, we immediately have $ h_{s,x}^{\varepsilon,\eta}(t)\lesssim \eta^{-\gamma_\varepsilon}<+\infty $. W.l.o.g., we take $s=0$ for simplicity and assume $T$ is \textit{small}, in particular $T\le 1 $. Let us write for $0<t\le T,x,y,y'\in \R^d $ the following:
	\begin{align*}
		&\Gamma (0,x,t,y')-\Gamma(0,x,t,y) = p_\alpha (t,y'-x)-p_\alpha(t,y-x)\\
		&\qquad + \int_{0}^{t/2}\int \Gamma (0,x,s,z)b(s,z)\cdot \left(\nabla_y p_\alpha (t-s,y-z)-\nabla_y p_\alpha (t-s,y'-z)\right) \d z \d s\\
		&\qquad + \int_{t/2}^{t-{(|y'-y|\vee\eta)^\alpha}}\int  \Gamma (0,x,s,z)b(s,z)\cdot \left(\nabla_y p_\alpha (t-s,y-z)-\nabla_{y'} p_\alpha (t-s,y'-z)\right) \d z \d s \\
		&\qquad + \int_{t-{(|y'-y|\vee\eta)^\alpha}}^t\int  \Gamma (0,x,s,z)b(s,z)\cdot \left(\nabla_y p_\alpha (t-s,y-z)-\nabla_{y'} p_\alpha (t-s,y'-z)\right) \d z \d s \\
		&=: \Delta_1 + \Delta_2 + \Delta_3+\Delta_4.
	\end{align*}
	{We {tacitly assume} as well that ${(|y'-y|\vee\eta)^\alpha}\le  t/2  $ since otherwise, i.e. in the off-diagonal case, the expected control $[(\Gamma (0,x,t,y)-\Gamma(0,x,t,y'))t^{\frac{\gamma_\varepsilon}\alpha}]/ [(\bar p_\alpha(t,y-x)+\bar p_\alpha(t,y'-x)) {(|y-y'| \vee \eta)^{\gamma_\eps}}]\le C$ readily follows from the Aronson type bounds \eqref{hold-ineq-density-diff}}.\\
	
	For $\Delta_1$, we use the regularity of the stable kernel, \eqref{hold-holder-space-palpha} to write
	\begin{equation}
		|\Delta_1|\lesssim \frac{|y-y'|^{\gamma_\varepsilon}}{t^{\frac{\gamma_\varepsilon}{\alpha}}}(\bar p_\alpha (t,y-x)+\bar p_\alpha(t,y'-x)).\\ \label{hold-maj-D1-holder-space-gamma}
	\end{equation}
	
	For $\Delta_2$, using again \eqref{hold-holder-space-palpha}, we write
	$$|\nabla_y p_\alpha (t-s,y-z)-\nabla_{y'} p_\alpha(t-s,y'-z)|\lesssim \frac{|y-y'|^{\gamma_\varepsilon}}{(t-s)^{\frac{\gamma_\varepsilon+1}{\alpha}}}(\bar p_\alpha (t-s,y-z)+\bar p_\alpha(t-s,y'-z))$$
	which yields, along with \eqref{hold-ineq-density-scheme} {and} the convolution property {\eqref{hold-APPR_CONV_PROP}} of $\bar p_\alpha$,
	\begin{align}
		|\Delta_2| &\lesssim \int_{0}^{t/2}\int \bar p_\alpha (s,z-x)\Vert b \Vert_{L^\infty} \frac{|y-y'|^{\gamma_\varepsilon}}{(t-s)^{\frac{\gamma_\varepsilon+1}{\alpha}}}(\bar p_\alpha (t-s,y-z)+\bar p_\alpha(t-s,y'-z))\d z \d s\nonumber\\
		&\lesssim {\frac{|y-y'|^{\gamma_\varepsilon}}{t^{\frac {\beta_\varepsilon}\alpha}}} (\bar p_\alpha (t,y-x)+\bar p_\alpha(t,y'-x)),\label{hold-maj-D2-holder-space-gamma}
	\end{align}
	where $\beta_\varepsilon:=\beta-\varepsilon$, noting that for $\gamma={\alpha+\beta-1}\in [1,2) $, ${\gamma_\varepsilon}+1-\alpha=2-(\alpha+\varepsilon)\le \beta-\varepsilon$ and recalling that we have assumed $T$ to be \textit{small}.\\
	
	For $\Delta_3$, using a Taylor expansion and then a cancellation argument, we have
	\begin{align}
		\Delta_3 &=  \int_{t/2}^{t-{(|y-y'|\vee\eta)^\alpha}}\int \int_0^1  \Gamma (0,x,s,z)b(s,z)\cdot \nabla_y^2 p_\alpha (t-s,y+\lambda (y'-y)-z)(y'-y)\d \lambda \d z \d s \notag\\
		&=  \int_{t/2}^{t-{(|y-y'|\vee\eta)^\alpha}}\int_0^1\int   [\Gamma (0,x,s,z)b(s,z)-\Gamma (0,x,s,y+\lambda (y'-y))b(s,y+\lambda (y'-y))]\notag\\ & \qquad \qquad\qquad\qquad \qquad \qquad \cdot \nabla_y^2 p_\alpha (t-s,y+\lambda (y'-y)-z)(y'-y) \d z \d \lambda \d s. \label{hold-DELTA_3_FOR_DIFF}
	\end{align}
	We then write
	\begin{align}
		&|\Gamma (0,x,s,z)b(s,z)-\Gamma (0,x,s,y+\lambda (y'-y))b(s,y+\lambda (y'-y))|\nonumber\\
		& \qquad \lesssim |y+\lambda (y'-y)-z|^\beta \bar p_\alpha(s,z-x) + \Vert b \Vert_{L^\infty}{\Big(h_{0,x}^{\varepsilon,\eta}(s) \frac{|y+\lambda (y'-y)-z|^{\gamma_\varepsilon}}{s^{\frac{{\gamma_\varepsilon}}{\alpha}}}{\mathbb 1}_{|y+\lambda (y'-y)-z|\le s^{\frac 1\alpha} }}\notag\\
		&\qquad\qquad {+\frac{|y+\lambda (y'-y)-z|^{\beta}}{s^{\frac{\beta}{\alpha}}}{\mathbb 1}_{|y+\lambda (y'-y)-z|\ge s^{\frac 1\alpha} }\Big) \left(\bar p_\alpha (s,z-x)+\bar p_\alpha (s,y+\lambda (y'-y)-x)\right)}\nonumber\\
		& \qquad \lesssim {|y+\lambda (y'-y)-z|^\beta \left(1+{s^{-\frac \beta\alpha}\left(1+\sup_{r\in (0,T]}h_{0,x}^{\varepsilon,\eta}(r)\right)} \right)\Bigg(\bar p_\alpha(s,z-x)   +\bar p_\alpha (s,y+\lambda (y'-y)-x) \Bigg)},\label{hold-pivot-post-cancel}
	\end{align}
	{recalling that $\beta \le \gamma_\varepsilon $ for the last inequality}.  	 Plugging this into {\eqref{hold-DELTA_3_FOR_DIFF}} and using \eqref{hold-drift-smoothing-iso-noise}, {recalling that on the considered time integration, $(t-s)\ge |y'-y|^\alpha$ (local diagonal regime)}, we get, 
	\begin{align*}
		|\Delta_3| &\lesssim \int_{t/2}^{t-{(|y-y'|\vee\eta)^\alpha}}\int_0^1\int \big(\bar p_\alpha(s,z-x)+\bar p_\alpha(s,y+\lambda(y'-y)-x)\big)  \left(1+{s^{-\frac \beta\alpha}\left(1+\sup_{r\in (0,T]}h_{0,x}^{\varepsilon,\eta}(r)\right)} \right)\\
		&\qquad \times \frac{|y-y'|^{\gamma_\varepsilon} }{(t-s)^{\frac{1+\gamma_\varepsilon-\beta}{\alpha}}}\bar p_\alpha (t-s,y+\lambda (y'-y)-z)\d z \d \lambda \d s \\
		&\lesssim \int_{t/2}^{t-{(|y-y'|\vee\eta)^\alpha}}\int_0^1 \bar p_\alpha(t,y+\lambda(y'-y)-x)  \left(1+{s^{-\frac \beta\alpha}\left(1+\sup_{r\in (0,T]}h_{0,x}^{\varepsilon,\eta}(r)\right)}  \right)\frac{|y-y'|^{\gamma_\varepsilon} }{(t-s)^{\frac{1+\gamma_\varepsilon-\beta}{\alpha}}}\d \lambda \d s\\
		&\lesssim \big(\bar p_\alpha(t,y-x)+\bar p_\alpha(t,y'-x)\big)\int_{t/2}^{t-{(|y-y'|\vee\eta)^\alpha}}   \left(1+{s^{-\frac \beta\alpha}\left(1+\sup_{r\in (0,T]}h_{0,x}^{\varepsilon,\eta}(r)\right)}  \right)\frac{|y-y'|^{\gamma_\varepsilon} }{(t-s)^{\frac{1+\gamma_\varepsilon-\beta}{\alpha}}} \d s,		  
	\end{align*}
	where, for the two last inequalities, we use the fact that for $s\geq t/2$, {up to a modification of the underlying variance in the Brownian case}, $\bar p_\alpha(s,y+\lambda(y'-y)-x)\lesssim \bar p_\alpha(t,y+\lambda(y'-y)-x)$ and since $|y-y'|\le t^{\frac 1\alpha} $, $\bar p_\alpha(t,y+\lambda(y'-y)-x)\le \bar p_\alpha(t,y-x)+\bar p_\alpha(t,y'-x) $ with the same previous abuse of notation if $\alpha=2 $. Finally, noting from the above definition of  $\gamma_\varepsilon {=(1\wedge \gamma)-\varepsilon}$ that $(1+\gamma_\varepsilon-\beta)/\alpha<1 {\iff \gamma>\gamma_\varepsilon}$, this yields
	\begin{align}
		|\Delta_3|
		&\lesssim \left(\bar p_\alpha(t,y-x) +\bar p_\alpha(t,y'-x) \right)
		|y-y'|^{\gamma_\eps}t^{\frac{\gamma-\gamma_\varepsilon}{\alpha}}\left(1+{t^{-\frac \beta\alpha}\left(1+\sup_{r\in (0,T]}h_{0,x}^{\varepsilon,\eta}(r)\right)}  \right).
		\label{hold-maj-D3-holder-space-gamma}
	\end{align}
	For $\Delta_4$, write:
	\begin{align*}
		|\Delta_4|\le &\left|\int_{t-{(|y'-y|\vee \eta)^\alpha}}^t\int  [\Gamma (0,x,s,z)b(s,z)- \Gamma (0,x,s,y)b(s,y)]\cdot \nabla_y p_\alpha (t-s,y-z)\d z \d s\right|\\
		&+\left|\int_{t-{(|y'-y|\vee \eta)^\alpha}}^t\int  [\Gamma (0,x,s,z)b(s,z)- \Gamma (0,x,s,y')b(s,y')]\cdot\nabla_{y'} p_\alpha (t-s,y'-z) \d z \d s\right|\\
		\lesssim& \int_{t-{(|y'-y|\vee \eta)^\alpha}}^{t}\int \bar p_{\alpha}(s,z-x)\left(1+{s^{-\frac \beta\alpha}\left(1+\sup_{r\in (0,T]}h_{0,x}^{\varepsilon,\eta}(r)\right)}  \right)\frac{\bar p_\alpha(t-s,y-z)+\bar p_\alpha(t-s,y'-z)}{(t-s)^{\frac 1\alpha-\frac\beta\alpha}} \d z\d s,
	\end{align*}
	where we used \eqref{hold-pivot-post-cancel} (with respectively $\lambda=0 $ and $\lambda=1 $ therein) {and \eqref{hold-drift-smoothing-iso-noise}} for the second inequality. We get:
	\begin{align}
		|\Delta_4|&{\lesssim}  (\bar p_\alpha(t,y-x)+\bar p_\alpha(t,y'-x)) {(|y'-y|\vee \eta)^\gamma}\left(1+{t^{-\frac \beta\alpha}\left(1+\sup_{r\in (0,T]}h_{0,x}^{\varepsilon,\eta}(r)\right)}  \right)\notag\\
		&{\lesssim}   (\bar p_\alpha(t,y-x)+\bar p_\alpha(t,y'-x)) {(|y'-y|\vee \eta)^{\gamma_\eps}}t^{\frac{\gamma-\gamma_\eps}{\alpha}}\left(1+{t^{-\frac \beta\alpha}\left(1+\sup_{r\in (0,T]}h_{0,x}^{\varepsilon,\eta}(r)\right)}  \right), \label{hold-maj-D4-holder-space-gamma}
	\end{align}	
	{where we also used the fact that {$(|y-y'| \vee \eta)^\alpha \leq t/2$} as previously mentioned for the last inequality}. Gathering estimates \eqref{hold-maj-D1-holder-space-gamma}, \eqref{hold-maj-D2-holder-space-gamma}, \eqref{hold-maj-D3-holder-space-gamma} and \eqref{hold-maj-D4-holder-space-gamma} {and considering the fact that {$|y-y'|^{\gamma_\eps}\leq (|y'-y|\vee \eta)^{\gamma_\eps}$}}, we obtain
	\begin{align*}
		|\Gamma (0,x,t,y)-\Gamma(0,x,t,y')| \lesssim &\left(\bar p_\alpha (t,y-x)+\bar p_\alpha (t,y'-x)\right)\frac{{(|y'-y|\vee \eta)^{\gamma_\eps}}}{t^{\frac{\gamma_\varepsilon}{\alpha}}}\\ & \quad \times \left[1 +t^{\frac{\gamma_\varepsilon-\beta_\varepsilon}{\alpha}}+t^{\frac{\gamma}{\alpha}}(1+t^{-\frac \beta\alpha})+t^{\frac{\gamma-\beta}{\alpha}}\sup_{s\in (0,T]}h_{0,x}^{\varepsilon,{\eta}}(s)\right].
	\end{align*}
	Noting that all the {exponents} of $t$ appearing {in}  brackets in the above equation are positive, we get, in turn taking the supremum for $t\in (0,T]$ on the l.h.s.
	\begin{align*}
		h_{0,x}^{\eps,{\eta}}(t) \lesssim 
		\left[1+T^{\frac{\gamma-\beta}{\alpha}}\sup_{r\in (0,T]}h_{0,x}^{\eps,{\eta}}(r)\right].
	\end{align*}
	Provided $T$ is small enough, {we obtain the (uniform in $\eta$) boundedness of $h_{0,x}^{\eps,{\eta}}$. Taking the limit $\eta \rightarrow 0$ concludes the proof of \eqref{hold-holder-space-gamma}.}\hfill $\square $

	\subsubsection{Proof of \eqref{hold-holder-space-gammah}: forward spatial H\"older regularity for the scheme}	

		{Let $\eps $ such that $\gamma_\varepsilon=(\gamma\wedge 1)-\varepsilon\ge\beta $ and set}
		$${g_{s,x}^{h,\varepsilon}(t)}:=\sup_{(z,z')\in (\R^d)^2} \left\{ \frac{|\Gamma^{h}(s,x,t,z)-\Gamma^{h}(s,x,t,z')|(t-s)^{\frac{\gamma_\varepsilon}{\alpha}}}{\left(p_\alpha(t-s,z-x)+p_\alpha(t-s,z'-x)\right)(|z-z'|{+h^{\frac 1\alpha}})^{\gamma_\varepsilon}} \right\}{.}$$
		{We emphasize that the time shift in the Duhamel representation of the scheme \eqref{hold-duhamel-scheme}, associated with the term $b(U_{\tau_s^h/h},X_{\tau_s}^h)\nabla_y p_\alpha(t-s,y-X_s^h) $, induces the additional term in ${h^{1/\alpha}}$ in the normalization. Intuitively, this can be explained since if $|y'-y|^\alpha\le h $, then, close to the time-boundary, i.e. for $s $ close to $t$, the local drift transition of the scheme of order $s-\tau_s^h $ is not negligible w.r.t $t-s$. {When looking at the diffusion, }this is usually dealt {with by} introducing a cut-off level {at $t-|y'-y|^\alpha $}. But on the remaining time interval, {one can still have $s-\tau_s^h\ge |y'-y|^\alpha $ and the drift somehow prevails for the scheme}}.\\
		
		Let us importantly {point out} that, from the Aronson type inequality \eqref{hold-ineq-density-scheme} for the scheme, it readily follows that ${g_{s,x}^{h,\varepsilon}(t)}\lesssim h^{-\frac{\gamma_\varepsilon}\alpha}<+\infty $. In particular, this means that this quantity can be used in a Gronwall or circular type procedure as we actually do below.\\
		
		Let us then introduce the following lemma:
		\begin{lemma}[Smoothing effect of the drift]\label{hold-THE_LEMMA_FOR_SCHEME} Let $\zeta$ be a multi-index with length $1\leq |\zeta|\leq {2}$. Then, for all $(x,\mathfrak{y})\in (\R^d)^2$, $0\leq t/2 \leq s {<} t \leq T$, $r>0$,
			\begin{align}
				&\left|\int \Gamma^{h}(0,x,\tau_s^h,z)b(r,z) \nabla_{\mathfrak{y}}^\zeta p_\alpha (t-{\tau_s^h},\mathfrak{y}-z) \d z\right|\nonumber\\
				& \qquad\qquad\qquad 
				\lesssim {\Big((t-\tau_s^h)^{\frac{\beta-|\zeta|}{\alpha}}+h^{\frac\beta\alpha} (t-\tau_s^h)^{-\frac{|\zeta|}{\alpha}}\Big)}\left(1+{\frac{1}{(\tau_s^h)^{\frac \beta\alpha}}}+\frac{g_{0,x}^{h,\varepsilon}(\tau_s^h)}{{(\tau_s^h)}^{\frac{\beta}{\alpha}}} \right) \bar p_\alpha (t,{\mathfrak{y}}-x).\label{hold-lemma-time-gain}
			\end{align}
		\end{lemma} 
		\begin{proof}[Proof of Lemma \ref{hold-THE_LEMMA_FOR_SCHEME}]
			To prove this, let us use the following cancellation argument
			\begin{align}
				I:=&\int \Gamma^{h}(0,x,\tau_s^h,z)b(r,z)\nabla_{\mathfrak{y}} p_\alpha (t-{\tau_s^h},\mathfrak{y}-z) \d z\nonumber\\
				& =  \int [\Gamma^{h}(0,x,\tau_s^h,z)b(r,z)-\Gamma^{h}(0,x,\tau_s^h,\mathfrak{y})b(r,\mathfrak{y})]\nabla_{\mathfrak{y}} p_\alpha (t-s,\mathfrak{y}-z) \d z{.}\label{hold-cancel-intra-lemma-hm}
			\end{align}
			Then, {we write similarly to \eqref{hold-pivot-post-cancel}},
			\begin{align}
				&|\Gamma^{h}(0,x,\tau_s^h,z)b(s,z)-\Gamma^{h}(0,x,\tau_s^h,{\mathfrak{y}})b(s,{\mathfrak{y}})|\nonumber\\
				&\qquad \lesssim |{\mathfrak{y}}-z|^\beta \bar p_\alpha(\tau_s^h,z-x) + \Vert b \Vert_{L^\infty}\Big({1+}g_{0,x}^{h,\varepsilon}(\tau_s^h)\Big) \frac{(|{\mathfrak{y}}-z|+h^{\frac 1\alpha})^{\beta}}{(\tau_s^h)^{\frac{\beta}{\alpha}}} \left(\bar p_\alpha (\tau_s^h,z-x)+\bar p_\alpha (\tau_s^h,{\mathfrak{y}}-x)\right)\nonumber\\
				&\qquad \lesssim  {(|{\mathfrak{y}}-z|+h^{\frac 1\alpha})^\beta  \left(1+\frac{1}{(\tau_s^h)^{\frac \beta\alpha}}+\frac{g_{0,x}^{h,\varepsilon}(\tau_s^h)}{{(\tau_s^h)}^{\frac{\beta}{\alpha}}} \right)\left(\bar p_\alpha(\tau_s^h,z-x)+\bar p_\alpha (\tau_s^h,{\mathfrak{y}}-x) \right).}\label{hold-pivot-post-cancel-4}
			\end{align}
			Plugging this into \eqref{hold-cancel-intra-lemma-hm} and using \eqref{hold-drift-smoothing-iso-noise}, we get 
			\begin{align*}
				&\left|\int \Gamma^{h}(0,x,\tau_s^h,z)b(r,z) \nabla_{\mathfrak{y}}^\zeta p_\alpha (t-{\tau_s^h},\mathfrak{y}-z) \d z\right|\nonumber\\
				\lesssim 	&{ \Big( (t-\tau_s^h)^{\frac{\beta-|\zeta|}{\alpha}}+h^{\frac \beta\alpha}(t-\tau_s^h)^{-\frac{|\zeta|}{\alpha}}\Big)}\left(1+{\frac{1}{(\tau_s^h)^{\frac \beta\alpha}}}+\frac{g_{0,x}^{h,\varepsilon}(\tau_s^h)}{{(\tau_s^h)}^{\frac{\beta}{\alpha}}} \right) \left[\bar p_\alpha (t,\mathfrak{y}-x) + \bar p_\alpha (\tau_s^h,{\mathfrak{y}}-x) \right].
			\end{align*}
			Using the fact that $s\geq t/2$,  {up to a modification of the underlying variance in the Brownian case}, $\bar p_\alpha(\tau_s^h,\mathfrak{y}-x)\lesssim \bar  p_\alpha(t,\mathfrak{y}-x)$, we obtain \eqref{hold-lemma-time-gain}.
		\end{proof}
		
		Let us first write
		\begin{align}
			&\Gamma^{h}(0,x,t,y')-\Gamma^{h}(0,x,t,y)  =p_\alpha (t,y'-x)-p_\alpha(t,y-x)\notag\\
			&\qquad + \int_{0}^{t/2}\frac{1}{h}\int_{\tau_s^h}^{\tau_s^h+h}\int \Gamma^{h} (0,x,\tau_s^h,z)b(r,z)\notag\\
			& \qquad \qquad \qquad \times{\E_{\tau_s^h,z,r}} \left[\nabla_y p_\alpha (t-s,y-X_s^{h})-\nabla_{y'} p_\alpha (t-s,y'-X_s^{h})\right]\d z \d r \d s\notag\\
			&\qquad + \int_{t/2}^{t}\frac{1}{h}\int_{\tau_s^h}^{\tau_s^h+h}\int \Gamma^{h} (0,x,\tau_s^h,z)b(r,z)\notag\\
			& \qquad \qquad \qquad \times {\E_{\tau_s^h,z,r}} \left[\nabla_y p_\alpha (t-s,y-X_s^{h})-\nabla_{y'} p_\alpha (t-s,y'-X_s^{h})\right] \d z \d r \d s\notag\\
			&=: \Delta_1 + \Delta_2 + \Delta_3,\label{hold-THE_DECOUP_HOLD_SCHEME}
		\end{align}
		{where we denoted $\E_{\tau_s^h,z,r}[\cdot]:=\E[\cdot|X_{\tau_s^h}^h=z,U_{\tau_s^h/h}=r] $}. \color{black} For $\Delta_1$, we use \eqref{hold-holder-space-palpha} to write
		\begin{align}
			{|\Delta_1|}=|p_\alpha (t,y'-x)-p_\alpha(t,y-x)|&\nonumber\lesssim \frac{|y-y'|^{{\gamma_\varepsilon}}}{t^{\frac{{\gamma_\varepsilon}}{\alpha}}}(\bar p_\alpha (t,y-x)+\bar p_\alpha(t,y'-x))\\&\lesssim \frac{(|y-y'|^\alpha\vee h)^{\frac{\gamma_\eps}{\alpha}}}{t^{\frac{{\gamma_\varepsilon}}{\alpha}}}(\bar p_\alpha (t,y-x)+\bar p_\alpha(t,y'-x)){.}
			\label{hold-maj-delta-1-holder-gammah}
		\end{align}
		For $\Delta_2$, we use \eqref{hold-holder-space-palpha} to write
		$$|\nabla_y p_\alpha (t-s,y-w)-\nabla_{y'} p_\alpha(t-s,y'-w)|\lesssim \frac{|y-y'|^{{\gamma_\varepsilon}}}{(t-s)^{\frac{{\gamma_\varepsilon}+1}{\alpha}}}(\bar p_\alpha (t-s,y-w)+\bar p_\alpha(t-s,y'-w))$$
		which yields, {similarly to \eqref{hold-maj-D2-holder-space-gamma}},
		\begin{align}
			|\Delta_2| & \lesssim  \int_{0}^{t/2}\int \bar p_\alpha (\tau_s^h,z-x) \Vert b \Vert_{L^\infty}\notag\\
			& \qquad \qquad \qquad \times \int \bar p_\alpha (s-\tau_s^h,w-z) \frac{|y-y'|^{\gamma_\varepsilon}}{(t-s)^{\frac{\gamma_\varepsilon+1}{\alpha}}}(\bar p_\alpha (t-s,y-w)+\bar p_\alpha(t-s,y'-w))\d w\d z \d s\notag\\
			&\lesssim \frac{(|y-y'|^\alpha\vee h)^{\frac{\gamma_\eps}{\alpha}}}{t^{\frac{\beta_\varepsilon}\alpha}}(\bar p_\alpha (t,y-x)+\bar p_\alpha(t,y'-x)),\label{hold-maj-delta-2-holder-gammah}
		\end{align}
		using the fact that the terms in $(t-s)$ can be taken out from the integral on the considered time interval.\\
		
		\color{black}
		For $\Delta_3$, let us develop the conditional expectation as follows: 
		\begin{align}
			&\E_{\tau_s^h,z,r} \left[\nabla_y p_\alpha (t-s,y-X_s^{h})-\nabla_{y'} p_\alpha (t-s,y'-X_s^{h})\right] \notag\\
			&= \E_{\tau_s^h,z,r} \left[\nabla_y p_\alpha \left(t-s,y-(z+b\left(r,z\right)(s-\tau_s^h)+Z_s-Z_{\tau_s^h})\right)\right]\notag\\
			& \quad -\E_{\tau_s^h,z,r} \left[\nabla_{y'} p_\alpha \left(t-s,y'-(z+b\left(r,z\right)(s-\tau_s^h)+Z_s-Z_{\tau_s^h})\right)\right]\notag\\ 
			&= \nabla_y p_\alpha \left(t-\tau_s^h,y-z-b\left(r,z\right)(s-\tau_s^h)\right)-\nabla_{y'} p_\alpha \left(t-\tau_s^h,y'-z-b\left(r,z\right)(s-\tau_s^h)\right)\notag\\ 
			&= \nabla_{w_2} p_\alpha \left(t-\tau_s^h,w_2\right)|_{w_2=y-z-b\left(r,z\right)(s-\tau_s^h)}-\nabla_{w_1} p_\alpha \left(t-\tau_s^h,w_1\right)|_{w_1=y-z-b\left(r,y\right)(s-\tau_s^h)}\notag\\
			&\quad + \nabla_{w_1} p_\alpha \left(t-\tau_s^h,w_1\right)|_{w_1=y-z-b\left(r,y\right)(s-\tau_s^h)}-\nabla_{w_1'} p_\alpha \left(t-\tau_s^h,w_1'\right)|_{w_1'=y'-z-b\left(r,y'\right)(s-\tau_s^h)}\notag\\
			&\quad + \nabla_{w_1'} p_\alpha \left(t-\tau_s^h,w_1'\right)|_{w_1'=y'-z-b\left(r,y'\right)(s-\tau_s^h)}-\nabla_{w_2'} p_\alpha \left(t-\tau_s^h,w_2'\right)|_{w_2'=y'-z-b\left(r,z\right)(s-\tau_s^h)}\notag\\
			&=-\int_0^1\nabla_{w_1}^2 p_\alpha \left(t-\tau_s^h,w_1-\mu [b(r,z)-b(r,y)] (s-\tau_s^h)\right)|_{w_1=y-z-b(r,y)(s-\tau_s^h)} \cdot\left[b\left(r,z\right)-b\left(r,y\right)\right](s-\tau_s^h)  \d \mu\notag\\
			& \quad+ \nabla_{w_1} p_\alpha \left(t-\tau_s^h,w_1\right)|_{w_1=y-z-b\left(r,y\right)(s-\tau_s^h)}-\nabla_{w_1'} p_\alpha \left(t-\tau_s^h,w_1'\right)|_{w_1'=y'-z-b\left(r,y'\right)(s-\tau_s^h)}\notag\\
			&\quad+ \int_0^1\nabla_{w_1'}^2 p_\alpha \left(t-\tau_s^h,w_1'-\mu [b(r,z)-b(r,y')] (s-\tau_s^h)\right)|_{w_1'=y'-z-b(r,y')(s-\tau_s^h)} \cdot\left[b\left(r,z\right)-b\left(r,y'\right)\right](s-\tau_s^h)  \d \mu,
			\label{hold-THE_DECOUP_DELTA_3_SCHEME}
		\end{align}
		yielding the corresponding terms $\Delta_{31}, \Delta_{32} $ and $\Delta_{33} $ once plugged into \eqref{hold-THE_DECOUP_HOLD_SCHEME}.
		For $\Delta_{31}$ and $\Delta_{33}$, since $t-\tau_s^h \geq s-\tau_s^h$, using the boundedness and the H\"older regularity of $b$ and then \eqref{hold-drift-smoothing-iso-noise-diag},  for $\tilde y\in \{y,y'\} $,
		\begin{align}
			&|\nabla_{\tilde w}^2 p_\alpha (t-\tau_s^h,\tilde w-\mu \left[b(r,z)-b(r,\tilde y)\right](s-\tau_s^h))|_{\tilde w=\tilde y-z-b(r,\tilde y)(s-\tau_s^h)}  	\cdot\left[b\left(r,z\right)-b\left(r,\tilde y\right)\right]  \notag\\
			\lesssim& (t-\tau_s^h)^{\frac{\beta-2}{\alpha}}\bar p_\alpha(t-\tau_s^h,\tilde y-z).\label{hold-KILLING_THE_DRIFT}
		\end{align}
		This then yields, using again $s-\tau_s^h\leq t-\tau_s^h$,
		\begin{align*}
			|\Delta_{31}| + |\Delta_{33}| & \lesssim\int_{t/2}^{t}\int \bar{p}_\alpha (\tau_s^h,z-x) (s-\tau_s^h)(t-\tau_s^h)^{\frac{\beta-2}{\alpha}}\left(\bar p_\alpha(t-\tau_s^h, y-z)+\bar p_\alpha(t-\tau_s^h, y'-z)\right)\d z  \d s\\
			& \lesssim  \left( \bar{p}_\alpha (t,y-x)+ \bar{p}_\alpha (t,y'-x)\right)h^{\frac{\gamma_\eps}{\alpha}}\int_{t/2}^{t}  (t-\tau_s^h)^{1+\frac{\beta-2-\gamma_\eps}{\alpha}}\d s.
		\end{align*}
		Note that, by definition, $\gamma_\eps =\big((\alpha + \beta-1) \wedge 1\big)-\eps \leq \alpha + \beta-1 -\eps$, so that
		$$1+\frac{\beta-2-\gamma_\eps}{\alpha}=\frac{\alpha+\beta-2-\gamma_\eps}{\alpha}=\frac{\gamma-\gamma_\varepsilon-1}{\alpha}\geq\frac{-1+\eps}{\alpha} >-1,$$
		using as well that $\alpha>1$ for the last inequality. In turn, we obtain
		\begin{align}
			|\Delta_{31}| + |\Delta_{33}|  &\lesssim  \left( \bar{p}_\alpha (t,y-x)+ \bar{p}_\alpha (t,y'-x)\right)h^{\frac{\gamma_\eps}{\alpha}}t^{\frac{\gamma-\gamma_\eps + \alpha-1}{\alpha}}\nonumber\\
			&\lesssim  \left( \bar{p}_\alpha (t,y-x)+ \bar{p}_\alpha (t,y'-x)\right)(|y-y'|^\alpha\vee h)^{\frac{\gamma_\eps}{\alpha}}T^{\frac{\gamma-\gamma_\eps + \alpha-1}{\alpha}}.\label{hold-maj-delta-3133-holder-gammah}
		\end{align}
		Let us turn to $\Delta_{32}$, which we split into two parts depending on whether the inner gradient is in diagonal (and in that case using a Taylor expansion) or off-diagonal regime:
		\begin{align*}
			\Delta_{32} =& \int_{t/2}^{t-|y-y'|^\alpha\vee h}\frac{1}{h}\int_{\tau_s^h}^{\tau_s^h+h}\int \Gamma^{h} (0,x,\tau_s^h,z)b(r,z)\\ & \times \int_{0}^{1}\nabla_{w_1}^2 p_\alpha \left(t-\tau_s^h,w_1+\mu \left[y'-y-(s-\tau_s^h)(b\left(r,y'\right)-b\left(r,y\right))\right]\right)|_{w_1=y-z-b\left(r,y\right)(s-\tau_s^h)} \\ & \qquad\qquad \cdot  \left[y'-y-(s-\tau_s^h)(b\left(r,y'\right)-b\left(r,y\right))\right]\d \mu  \d z \d r \d s\\
			& + \int_{t-|y-y'|^\alpha\vee h}^t\frac{1}{h}\int_{\tau_s^h}^{\tau_s^h+h}\int \Gamma^{h} (0,x,\tau_s^h,z)b(r,z)\\ &  \times\Bigg[\nabla_{w_1} p_\alpha \left(t-\tau_s^h,w_1\right)|_{w_1=y-z-b\left(r,y\right)(s-\tau_s^h)}-\nabla_{w_1'} p_\alpha \left(t-\tau_s^h,w_1'\right)_{w_1'=y-z-b\left(r,y\right)(s-\tau_s^h)}\Bigg]\d z \d r \d s\\
			=:&\Delta_{321}+\Delta_{322}.
		\end{align*}
		For $\Delta_{321}$, we use \eqref{hold-lemma-time-gain} with $\mathfrak{y}=y-b(r,y)(s-\tau_s^h)+\mu \left[y'-y-(b\left(r,y'\right)-b\left(r,y\right))(s-\tau_s^h)\right]$.
		Note that the term in $h^{\beta/\alpha}\le (t-\tau_s^h)^{\beta/\alpha} $ on the associated time regime for $\Delta_{321} $.			 
		
		Using as well the fact that, in the regime $|y-y'|\leq t^{\frac{1}{\alpha}}$, $$\bar{p}_\alpha(t,y-x-b\left(r,y\right)(s-\tau_s^h)+\mu \left[y'-y-(b\left(r,y'\right)-b\left(r,y\right))(s-\tau_s^h)\right])\lesssim \bar{p}_\alpha (t,y-x)$$ to obtain
		\begin{align*}
			|\Delta_{321}|&\lesssim  \bar{p}_\alpha (t,y-x) \int_{t/2}^{t-|y-y'|^\alpha\vee h} \frac{\left[|y-y'|+(s-\tau_s^h)|y-y'|^\beta\right]}{(t-\tau_s^h)^{\frac{2-\beta}{\alpha}}} \bar p_\alpha \left(t,y-x\right)\left(1+\frac{1}{(\tau_s^h)^{\frac \beta\alpha}}+\frac{g_{0,x}^{h,\varepsilon}(\tau_s^h)}{{(\tau_s^h)}^{\frac{\beta}{\alpha}}} \right)  \d s.
		\end{align*}
		%
		Recall that on the considered time interval $t-s\ge |y'-y|^\alpha\vee h $ and $|y-y'|\leq t^{\frac{1}{\alpha}} $. From	the fact that $s-\tau_s^h\leq t-\tau_s^h$,		
		we then get
		\begin{align}
			|\Delta_{321}|&\lesssim \bar p_\alpha \left(t,y-x\right)  (h^{\frac{\gamma_\eps}{\alpha}}t^{\frac{\beta}{\alpha}}+|y'-y|^{\gamma_\varepsilon})\int_{t/2}^{t-|y-y'|^\alpha\vee h} {(t-\tau_s^h)^{\frac{\beta-1-\gamma_\eps}{\alpha}}} \left(1+\frac{1}{(\tau_s^h)^{\frac \beta\alpha}}+\frac{g_{0,x}^{h,\varepsilon}(\tau_s^h)}{{(\tau_s^h)}^{\frac{\beta}{\alpha}}} \right) \d s\nonumber\\
			&\lesssim \bar p_\alpha \left(t,y-x\right)  (h^{\frac{\gamma_\eps}{\alpha}}+|y'-y|^{\gamma_\varepsilon})t^{\frac{\gamma-\gamma_\eps}{\alpha}}\left(1+t^{-\frac{\beta}{\alpha}}\left(1+\sup_{r\in (0,T]}g_{0,x}^{h,\varepsilon}(r)\right)\right),\label{hold-maj-delta-321-holder-gammah}
		\end{align}
		using as well $\frac{\beta-1-\gamma_\eps}{\alpha}>-1\iff \gamma>\gamma_\varepsilon $ for the last inequality.
		
		For $\Delta_{322}$, denote, for $\mathfrak{y}\in \{y,y'\}$,
		\begin{align*}
			&\delta_{322}(\mathfrak{y}):=\int_{t-|y-y'|^\alpha\vee h}^t\frac{1}{h}\int_{\tau_s^h}^{\tau_s^h+h} \\
			& \qquad \qquad \times \left| \int \Gamma^{h} (0,x,\tau_s^h,z)b(r,z) \cdot\E_{\tau_s^h,z} \left[\nabla_{\mathfrak{y}} p_\alpha \left(t-\tau_s^h,\mathfrak{y}-z-b\left(r,\mathfrak{y}\right)(s-\tau_s^h)\right)\right]\d z \right| \d r \d s
		\end{align*}
		so that $|\Delta_{322}|\leq \delta_{322}(y)+\delta_{322}(y')$. Similarly to $\Delta_{321}$, we now use \eqref{hold-lemma-time-gain} with $\mathfrak{y}=\mathfrak{y}-b\left(r,\mathfrak{y}\right)(s-\tau_s^h)$ observing that on the considered time integration interval the term in $h^{\beta/\alpha} $ of \eqref{hold-lemma-time-gain} remains:
		\begin{align}
			\delta_{322}(\mathfrak{y})&\lesssim\bar{p}_\alpha (t,\mathfrak{y}-x)\int_{t-|y-y'|^\alpha\vee h}^t\Big((t-\tau_s^h)^{\frac{\beta - 1}{\alpha}}+h^{\frac\beta\alpha}(t-\tau_s^h)^{-\frac 1\alpha}\Big)\left(1+\frac{1}{(\tau_s^h)^{\frac \beta\alpha}}+\frac{g_{0,x}^{h,\varepsilon}(\tau_s^h)}{{(\tau_s^h)}^{\frac{\beta}{\alpha}}} \right) \d s\nonumber\\
			&\lesssim\bar{p}_\alpha (t,\mathfrak{y}-x)(|y-y'|^\alpha \vee h)^{\frac{\gamma}{\alpha}}\left(1+t^{-\frac{\beta}{\alpha}}\left(1+\sup_{r\in (0,T]}g_{0,x}^{h,\varepsilon}(r)\right)\right).\label{hold-maj-delta-322-holder-gammah}
		\end{align}
		Gathering \eqref{hold-maj-delta-3133-holder-gammah}, \eqref{hold-maj-delta-321-holder-gammah} and \eqref{hold-maj-delta-322-holder-gammah}, we obtain
		\begin{align}
			|\Delta_{3}| &\lesssim \left(\bar p_\alpha(t,y-x)+\bar{p}_\alpha(t,y'-x)\right)(|y-y'|^\alpha \vee h)^{\frac{\gamma_\eps}{\alpha}}\nonumber\\ & \qquad \qquad \qquad \times \left(1+t^{-\frac{\beta}{\alpha}}\left(1+\sup_{r\in (0,T]}g_{0,x}^{h,\varepsilon}(r)\right)\right) \left[(|y-y'|^\alpha \vee h)^{\frac{\gamma-\gamma_\eps}{\alpha}}+t^{\frac{\gamma-\gamma_\eps}{\alpha}}+t^{\frac{\gamma-\gamma_\eps + \alpha-1}{\alpha}}\right].\label{hold-maj-delta-32-holder-gammah}
		\end{align}
		Along with \eqref{hold-maj-delta-1-holder-gammah} and \eqref{hold-maj-delta-2-holder-gammah}, we eventually get
		\begin{align*}
			& \frac{|\Gamma^{h}(0,x,t,y)-\Gamma^{h}(0,x,t,y')|t^{\frac{\gamma_\eps}{\alpha}}}{\left(\bar p_\alpha(t,y-x)+\bar{p}_\alpha(t,y'-x)\right)(|y-y'|^\alpha \vee h)^{\frac{\gamma_\eps}{\alpha}}}  \lesssim 1+T^{\frac{\gamma_\eps-\beta_\eps}{\alpha}} \nonumber \\ & \qquad\qquad \qquad  + \left(1+t^{-\frac{\beta}{\alpha}}\left(1+\sup_{r\in (0,T]}g_{0,x}^{h,\varepsilon}(r)\right)\right)  \left[(|y-y'|^\alpha \vee h)^{\frac{\gamma-\gamma_\eps}{\alpha}}t^{\frac{\gamma_\eps}{\alpha}}+t^{\frac{\gamma}{\alpha}}+t^{\frac{\gamma+ \alpha-1}{\alpha}}\right]\\
			&\qquad\qquad \qquad\lesssim1+T^{\frac{\gamma_\eps-\beta_\eps}{\alpha}} +T^{\frac{\gamma-\beta}\alpha}+ \sup_{r\in (0,T]}g_{0,x}^{h,\varepsilon}(r)  \left[t^{\frac{\gamma-\beta}{\alpha}}+t^{\frac{\gamma+ \alpha-1-\beta}{\alpha}}\right],
		\end{align*}
		using as well that  $ |y'-y|^\alpha\vee h \le t $ for the last inequality. Since $ \gamma-\beta=\alpha-1>0$, Equation \eqref{hold-holder-space-gammah} then follows taking the supremum in time in the previous inequality provided $T $ is small enough. \hfill $\square $
		\color{black}

		\subsection{Proof of \eqref{hold-holder-time-gamma}: forward time H\"older regularity for the diffusion}
		\begin{proof}[Proof of \eqref{hold-holder-time-gamma}] We proceed here with the proof of the forward time sensitivity. Importantly, the proof will use the previously proved claim \eqref{hold-holder-space-gamma} which actually also gives $\beta $-H\"older sensitivity since $\beta\le \gamma_\varepsilon $. 
			W.l.o.g. we take $s=0$ for notational simplicity. Starting from the Duhamel representation \eqref{hold-duhamel-Diff} and using cancellation arguments, we have {for $0< t<t'\le T,\ x,y\in \R^d$},
			\begin{align}
				\Gamma(0,x,t,y)-\Gamma(0,x,t',y)
				&= p_\alpha (t,y-x)-  p_\alpha(t',y-x)\notag\\
				&\qquad +\int_t^{t'} \int \Gamma(0,x,s,z)b(r,z)\cdot\nabla_y p_\alpha(t'-s,y-z)\d z \d s \nonumber\\
				&\qquad +\int_0^{t} \int \Gamma(0,x,s,z)b(r,z)\cdot[\nabla_y p_\alpha (t'-s,y-z)-\nabla_y p_\alpha (t-s,y-z)]\d z \d s \nonumber\\
				&= p_\alpha (t,y-x)-  p_\alpha(t',y-x)\notag\\
				&\qquad +\int_t^{t'} \int [\Gamma(0,x,s,z)b(r,z)-\Gamma(0,x,s,y)b(r,y)]\cdot\nabla_y p_\alpha(t'-s,y-z)\d z \d s \nonumber\\
				&\qquad +\int_0^{t} \int \Gamma(0,x,s,z)b(r,z)\nonumber\\
				& \qquad\qquad\qquad\qquad\qquad\qquad\cdot[\nabla_y p_\alpha (t'-s,y-z)-\nabla_y p_\alpha (t-s,y-z)]\d z \d s \nonumber\\
				&=:H_1 + H_2 + H_3.
			\end{align} 
			For $H_1$, we directly use \eqref{hold-holder-time-palpha} to write
			\begin{equation}
				|H_1|=|p_\alpha (t,y-x)-  p_\alpha(t',y-x)|\lesssim \frac{(t'-t)^{\frac{\gamma}{\alpha}}}{t^{\frac{\gamma}{\alpha}}} p_\alpha (t,y-x).\label{hold-maj-H1-holder-time-gamma}
			\end{equation}
			For $H_2$, let us use the regularity of $b$ and the forward spatial regularity of $\Gamma$, \eqref{hold-holder-space-gamma} to write:
			\begin{align}
				|\Gamma (0,x,s,z)b(s,z)-\Gamma(0,x,s,y)b(s,y)|	&\lesssim \bar p_\alpha(s,z-x) |y-z|^\beta  \left(1+\frac{1}{s^{\frac{\beta}{\alpha}}} \right)+\bar p_\alpha (s,y-x)  \frac{|y-z|^{\beta}}{s^{\frac{\beta}{\alpha}}}.\label{hold-pivot-post-cancel-2}
			\end{align}
			Plugging this into $H_2$ and using the fact that for $s\geq t$, $\bar p_\alpha (s,y-x)\lesssim \bar p_\alpha (t,y-x)$ along with \eqref{hold-drift-smoothing-iso-noise}, we have
			\begin{align}
				|H_2|&\lesssim \int_t^{t'} \int \bar p_\alpha(s,z-x) |y-z|^\beta  \left(1+\frac{1}{s^{\frac{\beta}{\alpha}}} \right) |\nabla_y p_\alpha(t'-s,y-z)|\d z \d s\nonumber\\
				& \qquad +  \int_t^{t'} \int \bar  p_\alpha (s,y-x)  \frac{|y-z|^{\beta}}{s^{\frac{\beta}{\alpha}}} |\nabla_y p_\alpha(t'-s,y-z)|\d z \d s\nonumber\\
				&\lesssim \int_t^{t'} \int \bar p_\alpha(s,z-x)  \left(1+\frac{1}{s^{\frac{\beta}{\alpha}}} \right)(t'-s)^{\frac{\beta-1}{\alpha}}{\bar p_\alpha (t'-s,y-z)}\d z \d s\nonumber\\
				& \qquad +  \bar p_\alpha (t,y-x)\int_t^{t'} \int   \frac{(t'-s)^{\frac{\beta-1}{\alpha}}}{s^{\frac{\beta}{\alpha}}}  {\bar p_\alpha}(t'-s,y-z)|\d z \d s\nonumber\\
				&\lesssim \bar p_\alpha (t,y-x){(1+t^{-\frac \beta \alpha})}(t-t')^{\frac{\gamma}{\alpha}}.\label{hold-maj-H2-holder-time-gamma}
			\end{align}
			For $H_3$, let us split it into three parts again, using a cancellation argument on two of them:
			\begin{align*}
				H_3 &= \int_0^{\frac{t-(t'-t)}{2}} \int\Gamma(0,x,s,z)b(r,z) \cdot [\nabla_y p_\alpha (t'-s,y-z)-\nabla_y p_\alpha (t-s,y-z)]\d z \d s \\
				&\qquad + \int_{\frac{t-(t'-t)}{2}}^{t-(t'-t)} \int \left[\Gamma(0,x,s,z)b(r,z)-\Gamma(0,x,s,y)b(r,y)\right]\cdot [\nabla_y p_\alpha (t'-s,y-z)-\nabla_y p_\alpha (t-s,y-z)]\d z \d s \\
				&\qquad + \int_{t-(t'-t)}^t \int \left[\Gamma(0,x,s,z)b(r,z)-\Gamma(0,x,s,y)b(r,y)\right]\cdot[\nabla_y p_\alpha (t'-s,y-z)-\nabla_y p_\alpha (t-s,y-z)]\d z \d s \\
				&=: H_{31}+H_{32}+H_{33}.
			\end{align*}
			For $H_{31}$, notice that for $s\leq (t-(t'-t))/{2}=t-t'/2$, $t'-s \geq 3t'/2-t$ and $t-s \geq t'/2$, {there will be no time singularities in $(t-s)$ or $(t'-s)$ to integrate}.  There is {thus} no need to use a cancellation argument to derive a smoothing effect. Using simply \eqref{hold-holder-time-palpha}, we get
			\begin{align}
				|H_{31}|& \lesssim \int_0^{\frac{t-(t'-t)}{2}} \int \bar p_\alpha (s,z-x) \Vert b \Vert_{L^\infty}  \frac{(t'-t)^{\frac{\gamma}{\alpha}}}{(t-s)^{\frac{\gamma+1}{\alpha}}}[\bar  p_\alpha (t'-s,y-z)+\bar  p_\alpha (t-s,y-z)]\d z \d s\nonumber \\
				&\lesssim (t'-t)^{\frac{\gamma}{\alpha}} {t^{-\frac\beta\alpha}}\left[\bar p_\alpha (t,y-x)+\bar p_\alpha (t',y-x)\right].\label{hold-maj-H31-holder-time-gamma}
			\end{align}
			For $H_{32}$, notice that on the considered time interval, $t'-s\geq 2(t'-t)$, so we are at the right time scale to use a Taylor expansion in time:
			\begin{align*}
				&H_{32}\\
				=& \int_{\frac{t-(t'-t)}{2}}^{t-(t'-t)} \int \left[\Gamma(0,x,s,z)b(r,z)-\Gamma(0,x,s,y)b(r,y)\right]\cdot \int_0^1 \partial_t\nabla_y p_\alpha (t-s+\lambda (t'-t),y-z) (t'-t) \d \lambda\d z \d s{.}
			\end{align*}
			Using \eqref{hold-pivot-post-cancel-2} and then \eqref{hold-drift-smoothing-iso-noise} along with the fact that for on the considered time interval, we have $t-s+\lambda (t'-t)\asymp t-s$ and  $ \bar{p}_\alpha (s,y-x)\lesssim \bar{p}_\alpha (t,y-x)$, we get 
			\begin{align}
				|H_{32}|&\lesssim(t'-t) \int_{\frac{t-(t'-t)}{2}}^{t-(t'-t)} \int \left[\bar p_\alpha (s,z-x)\left(1+s^{-\frac{\beta}{\alpha}}\right)+\bar p_\alpha(s,y-x)s^{-\frac{\beta}{\alpha}} \right]\nonumber\\ & \qquad\qquad\qquad\qquad \times \int_0^1 |y-z|^\beta| \partial_t\nabla_y p_\alpha (t-s+\lambda (t'-t),y-z) | \d \lambda\d z \d s\nonumber\\
				&\lesssim(t'-t) \int_{\frac{t-(t'-t)}{2}}^{t-(t'-t)} \int \left[\bar p_\alpha (s,z-x)\left(1+s^{-\frac{\beta}{\alpha}}\right)+\bar p_\alpha(t,y-x)s^{-\frac{\beta}{\alpha}} \right] (t-s)^{-1+\frac{\beta-1}{\alpha}} \bar p_\alpha (t-s,y-z) | \d z \d s\nonumber\\
				&\lesssim (t'-t)^{\frac{\gamma}{\alpha}}{(1+t^{-\frac \beta\alpha})} \bar{p}_\alpha(t,y-x),\label{hold-maj-H32-holder-time-gamma}
			\end{align}
			{	observing that $-1+(\beta-1)/\alpha<-1 $ for the last inequality}.\\

			{For $H_{33}$, we take advantage of the fact that we integrate on a time interval whose length  corresponds to the time difference. This means we can use the smoothing in time effect for each term of the difference (no need to expand in space the difference of the gradients). Namely,}
			\begin{align*}
				|H_{33}|\le& \left|\int_{t-(t'-t)}^t \int \left[\Gamma(0,x,s,z)b(r,z)-\Gamma(0,x,s,y)b(r,y)\right]\cdot \nabla_y p_\alpha (t'-s,y-z)\right|\\
				&+\left|\int_{t-(t'-t)}^t \int \left[\Gamma(0,x,s,z)b(r,z)-\Gamma(0,x,s,y)b(r,y)\right]\cdot \nabla_y p_\alpha (t-s,y-z)\d z \d s\right|\\
				\le &\int_{t-(t'-t)}^t \int \left[\bar p_\alpha (s,z-x)\left(1+s^{-\frac{\beta}{\alpha}}\right)+\bar p_\alpha(t,y-x)s^{-\frac{\beta}{\alpha}} \right] (t-s)^{\frac{\beta-1}{\alpha}}\\
				&\qquad \times 
				\big(\bar p_\alpha (t-s,y-z)+\bar p_\alpha (t'-s,y-z)\big) \d z \d s,
			\end{align*}
			using again \eqref{hold-pivot-post-cancel-2} for the last inequality. This eventually yields
			\begin{equation}
				|H_{33}| \lesssim (t'-t)^{\frac{\gamma}{\alpha}} {(1+t^{-\frac{\beta}{\alpha}})}(\bar p_\alpha (t,y-x)+\bar p_\alpha (t',y-x)){.}\label{hold-maj-H33-holder-time-gamma}
			\end{equation}
			Gathering estimates \eqref{hold-maj-H1-holder-time-gamma}, \eqref{hold-maj-H2-holder-time-gamma}, \eqref{hold-maj-H31-holder-time-gamma}, \eqref{hold-maj-H32-holder-time-gamma} and \eqref{hold-maj-H33-holder-time-gamma}, we obtain
			\begin{equation}
				|\Gamma (0,x,t,y)-\Gamma(0,x,t',y)| \lesssim \left(\frac{t'-t}{t}\right)^{\frac{\gamma}{\alpha}}\left[\bar p_\alpha (t,y-x)+\bar p_\alpha (t',y-x)\right],
			\end{equation}
			which precisely gives  \eqref{hold-holder-time-gamma} since we have assumed $s=0$.
		\end{proof}

\appendix		
		\section{About the full parabolic bootstrap in the forward variable for the diffusion}
		\label{hold-BOOT} 	
		The point of this section is to provide a proof of the full parabolic bootstrap for the diffusion in its forward variable in the case $\alpha+\beta-1<1 $. Indeed, when $\alpha+\beta-1\ge 1 $, it cannot be expected to have an exponent greater than $1$ and \eqref{hold-holder-space-gamma} is already sharp.\\
		
		Namely, we prove the following : there exists $C :=C({d},b,\alpha,T)$ s.t. for all $0\le s<t\le T,\ (x,y,w)\in (\R^d)^3$ s.t.  $|y-w|\le (t-s)^{\frac 1\alpha} $,
		\begin{align}
			\label{hold-holder-space-gamma_SHARP}
			|\Gamma(s,x,t,y)-\Gamma(s,x,t,w)|\le {C}\left(\frac{|y-w|}{(t-s)^{\frac1\alpha}} \right)^{{\gamma}}\bar p_\alpha(t-s,w-x).
		\end{align}
		
		The approach is very similar to the previous one to show \eqref{hold-holder-space-gamma} and we present here the result for the sake of completeness only as we do not make use of it. Equation \eqref{hold-holder-space-gamma} (taking $\varepsilon $ therein s.t. $\gamma_\varepsilon=\beta $) is enough for the proof of Theorem \ref{hold-thm-main}.\\
		
		\begin{proof} Set, for $\eta $ meant to be small, 
			$$h_{s,x}^{\eta}(t):=\sup_{(z,z')\in (\R^d)^2} \left\{ \frac{|\Gamma(s,x,t,z)-\Gamma(s,x,t,z')|(t-s)^{\frac{\gamma}{\alpha}}}{\left(\bar p_\alpha(t-s,z-x)+\bar p_\alpha(t-s,z'-x)\right){(|z-z'| \vee \eta)^{\gamma}}} \right\}.$$
			{Since we already know from \eqref{hold-ineq-density-diff} that $ \frac{|\Gamma(s,x,t,z)-\Gamma(s,x,t,z')|}{\bar p_\alpha(t-s,z-x)+\bar p_\alpha(t-s,z'-x)}<\infty$}, we immediately have $ h_{s,x}^{\eta}(t)\lesssim \eta^{-\gamma}<+\infty $. W.l.o.g., we take $s=0$ for simplicity and assume $T$ is \textit{small}, in particular $T\le 1 $. Let us write for $0<t\le T,(x,y,y')\in (\R^d)^3 $ the following:
			\begin{align*}
				&\Gamma (0,x,t,y')-\Gamma(0,x,t,y) = p_\alpha (t,y'-x)-p_\alpha(t,y-x)\\
				&\qquad + \int_{0}^{t/2}\int \Gamma (0,x,s,z)b(s,z)\cdot \left(\nabla_y p_\alpha (t-s,y-z)-\nabla_{y'} p_\alpha (t-s,y'-z)\right) \d z \d s\\
				&\qquad + \int_{t/2}^{t-{(|y'-y|\vee\eta)^\alpha}}\int  \Gamma (0,x,s,z)b(s,z)\cdot \left(\nabla_y p_\alpha (t-s,y-z)-\nabla_{y'} p_\alpha (t-s,y'-z)\right) \d z \d s \\
				&\qquad + \int_{t-{(|y'-y|\vee\eta)^\alpha}}^t\int  \Gamma (0,x,s,z)b(s,z)\cdot \left(\nabla_y p_\alpha (t-s,y-z)-\nabla_{y'} p_\alpha (t-s,y'-z)\right) \d z \d s \\
				&=: \Delta_1 + \Delta_2 + \Delta_3+\Delta_4.
			\end{align*}
			{We {tacitly assume} as well that ${(|y'-y|\vee\eta)^\alpha}\le  t/2  $ since otherwise, i.e. in the off-diagonal case, the expected control $[(\Gamma (0,x,t,y)-\Gamma(0,x,t,y'))t^{\frac{\gamma}\alpha}]/ [(\bar p_\alpha(t,y-x)+\bar p_\alpha(t,y'-x)) {(|y-y'| \vee \eta)^{\gamma}}]\le C$ readily follows from the Aronson type bounds \eqref{hold-ineq-density-diff}}.\\
			
			For $\Delta_1$, we use the regularity of the stable kernel, \eqref{hold-holder-space-palpha} to write
			\begin{equation}
				|\Delta_1|\lesssim \frac{|y-y'|^{\gamma}}{t^{\frac{\gamma}{\alpha}}}(\bar p_\alpha (t,y-x)+\bar p_\alpha(t,y'-x)).\\ \label{hold-maj-D1-holder-space-gamma_S}
			\end{equation}
			
			For $\Delta_2$, using again \eqref{hold-holder-space-palpha}, we write
			$$|\nabla_y p_\alpha (t-s,y-z)-\nabla_{y'} p_\alpha(t-s,y'-z)|\lesssim \frac{|y-y'|^{\gamma}}{(t-s)^{\frac{\gamma+1}{\alpha}}}(\bar p_\alpha (t-s,y-z)+\bar p_\alpha(t-s,y'-z))$$
			which yields, along with \eqref{hold-ineq-density-scheme} {and} the convolution property {\eqref{hold-APPR_CONV_PROP}} of $\bar p_\alpha$,
			\begin{align}
				|\Delta_2| &\lesssim \int_{0}^{t/2}\int \bar p_\alpha (s,z-x)\Vert b \Vert_{L^\infty} \frac{|y-y'|^{\gamma}}{(t-s)^{\frac{\gamma+1}{\alpha}}}(\bar p_\alpha (t-s,y-z)+\bar p_\alpha(t-s,y'-z))\d z \d s\nonumber\\
				&\lesssim {\frac{|y-y'|^{\gamma}}{t^{\frac {\beta}\alpha}}} (\bar p_\alpha (t,y-x)+\bar p_\alpha(t,y'-x)).\label{hold-maj-D2-holder-space-gamma_S}
			\end{align}
			
			For $\Delta_3$, using a Taylor expansion and then a cancellation argument, we have
			\begin{align}
				\Delta_3 &=  \int_{t/2}^{t-{(|y-y'|\vee\eta)^\alpha}}\int \int_0^1  \Gamma (0,x,s,z)b(s,z)\cdot \nabla_y^2 p_\alpha (t-s,y+\lambda (y'-y)-z)(y'-y)\d \lambda \d z \d s \notag\\
				&=  \int_{t/2}^{t-{(|y-y'|\vee\eta)^\alpha}}\int_0^1\int   [\Gamma (0,x,s,z)b(s,z)-\Gamma (0,x,s,y+\lambda (y'-y))b(s,y+\lambda (y'-y))]\notag\\ & \qquad \qquad\qquad\qquad \qquad \qquad \cdot \nabla_y^2 p_\alpha (t-s,y+\lambda (y'-y)-z)(y'-y) \d z \d \lambda \d s. \label{hold-DELTA_3_FOR_DIFF_S}
			\end{align}
			We then write, using \eqref{hold-holder-space-gamma} with $\varepsilon>0$ s.t. $\gamma_\eps=\beta $ when $|y+\lambda (y'-y)-z|\le s^{\frac 1\alpha}$ and a triangular inequality when $|y+\lambda (y'-y)-z|\ge s^{\frac 1\alpha} $,
			\begin{align}
				&|\Gamma (0,x,s,z)b(s,z)-\Gamma (0,x,s,y+\lambda (y'-y))b(s,y+\lambda (y'-y))|\nonumber\\
				& \qquad \lesssim |y+\lambda (y'-y)-z|^\beta \bar p_\alpha(s,z-x) \notag\\
				&\qquad\qquad + \Vert b \Vert_{L^\infty}{\frac{|y+\lambda (y'-y)-z|^{\beta}}{s^{\frac{{\beta}}{\alpha}}}{\mathbb 1}_{|y+\lambda (y'-y)-z|\le s^{\frac 1\alpha} }}\left(\bar p_\alpha (s,z-x)+\bar p_\alpha (s,y+\lambda (y'-y)-x)\right)\nonumber\\
				& \qquad \lesssim {|y+\lambda (y'-y)-z|^\beta \left(1+{s^{-\frac \beta\alpha}} \right)\Bigg(\bar p_\alpha(s,z-x)   +\bar p_\alpha (s,y+\lambda (y'-y)-x) \Bigg)}.\label{hold-pivot-post-cancel_S}
			\end{align} 
			Plugging this into {\eqref{hold-DELTA_3_FOR_DIFF_S}} and using \eqref{hold-drift-smoothing-iso-noise}, {recalling that on the considered time integration $(t-s)\ge |y'-y|^\alpha$ (local diagonal regime)}, {we get for $\gamma_1\in (\gamma,1) $},
			\begin{align*}
				|\Delta_3| &\lesssim \int_{t/2}^{t-{(|y-y'|\vee\eta)^\alpha}}\int_0^1\int \big(\bar p_\alpha(s,z-x)+\bar p_\alpha(s,y+\lambda(y'-y)-x)\big)  \left(1+{s^{-\frac \beta\alpha}} \right)\\
				&\qquad \times \frac{|y-y'|^{\gamma_1} }{(t-s)^{\frac{1+\gamma_1-\beta}{\alpha}}}\bar p_\alpha (t-s,y+\lambda (y'-y)-z)\d z \d \lambda \d s \\
				&\lesssim \int_{t/2}^{t-{(|y-y'|\vee\eta)^\alpha}}\int_0^1 \bar p_\alpha(t,y+\lambda(y'-y)-x)  \left(1+{s^{-\frac \beta\alpha}}\right)\frac{|y-y'|^{\gamma_1} }{(t-s)^{\frac{1+\gamma_1-\beta}{\alpha}}}\d \lambda \d s\\
				&\lesssim \big(\bar p_\alpha(t,y-x)+\bar p_\alpha(t,y'-x)\big)\int_{t/2}^{t-{(|y-y'|\vee\eta)^\alpha}}   \left(1+{s^{-\frac \beta\alpha}}  \right)\frac{|y-y'|^{\gamma_1} }{(t-s)^{\frac{1+\gamma_1-\beta}{\alpha}}} \d s,		  
			\end{align*}
			where, for the two last inequalities, we use the fact that for $s\geq t/2$, {up to a modification of the underlying variance in the Brownian case}, $\bar p_\alpha(s,y+\lambda(y'-y)-x)\lesssim \bar p_\alpha(t,y+\lambda(y'-y)-x)$ and since $|y-y'|\le t^{\frac 1\alpha} $, $\bar p_\alpha(t,y+\lambda(y'-y)-x)\le \bar p_\alpha(t,y-x)+\bar p_\alpha(t,y'-x) $ with the same previous abuse of notation if $\alpha=2 $. Finally, noting from the above $\gamma_1$,   that $(1+\gamma_1-\beta)/\alpha>1 {\iff \gamma_1>\gamma}$, this yields
			\begin{align}
				|\Delta_3|
				&\lesssim \left(\bar p_\alpha(t,y-x) +\bar p_\alpha(t,y'-x) \right)
				|y-y'|^{\gamma_1}{(|y-y'|\vee \eta)^{\gamma-\gamma_1}}\lesssim
				\left(\bar p_\alpha(t,y-x) +\bar p_\alpha(t,y'-x) \right)
				(|y-y'|\vee \eta)^{\gamma}.
				\label{hold-maj-D3-holder-space-gamma_S}
			\end{align}
			For $\Delta_4$, write:
			\begin{align*}
				|\Delta_4|\le &\left|\int_{t-{(|y'-y|\vee \eta)^\alpha}}^t\int  [\Gamma (0,x,s,z)b(s,z)- \Gamma (0,x,s,y)b(s,y)]\cdot \nabla_y p_\alpha (t-s,y-z)\d z \d s\right|\\
				&+\left|\int_{t-{(|y'-y|\vee \eta)^\alpha}}^t\int  [\Gamma (0,x,s,z)b(s,z)- \Gamma (0,x,s,y')b(s,y')]\cdot\nabla_{y'} p_\alpha (t-s,y'-z) \d z \d s\right|\\
				\lesssim& \int_{t-{(|y'-y|\vee \eta)^\alpha}}^{t}\int \bar p_{\alpha}(s,z-x)\left(1+{s^{-\frac \beta\alpha}}  \right)\frac{\bar p_\alpha(t-s,y-z)+\bar p_\alpha(t-s,y'-z)}{(t-s)^{\frac 1\alpha-\frac\beta\alpha}} \d s\d z,
			\end{align*}
			where we used \eqref{hold-holder-space-gamma}, with $\gamma_\varepsilon=\beta $ and \eqref{hold-drift-smoothing-iso-noise} for the second inequality. We get:
			\begin{align}
				|\Delta_4|&{\lesssim}  (\bar p_\alpha(t,y-x)+\bar p_\alpha(t,y'-x)) {(|y'-y|\vee \eta)^\gamma}\left(1+{t^{-\frac \beta\alpha}}\right).
				\label{hold-maj-D4-holder-space-gamma_S}
			\end{align}	
			Gathering estimates \eqref{hold-maj-D1-holder-space-gamma_S}, \eqref{hold-maj-D2-holder-space-gamma_S}, \eqref{hold-maj-D3-holder-space-gamma_S} and \eqref{hold-maj-D4-holder-space-gamma_S}, we obtain
			\begin{align*}
				|\Gamma (0,x,t,y)-\Gamma(0,x,t,y')| \lesssim &\left(\bar p_\alpha (t,y-x)+\bar p_\alpha (t,y'-x)\right)\frac{{(|y'-y|\vee \eta)^{\gamma}}}{t^{\frac{\gamma}{\alpha}}} \left[1 +t^{\frac{\gamma-\beta}{\alpha}}\right].
			\end{align*}
			Noting that $\gamma-\beta=\alpha-1>0 $, we get, in turn
			\begin{align*}
				h_{0,x}^{{\eta}}(t) \lesssim 
				1.
			\end{align*}
			Taking the limit $\eta \rightarrow 0$ concludes the proof of \eqref{hold-holder-space-gamma_SHARP}.
			
		\end{proof}

\section{Some technical results about Besov spaces}		
		
		 We first recall that \textcolor{black}{denoting by ${\mathcal S'}( \R^d) $ the dual space of the Schwartz class ${\mathcal S}( \R^d) $}, for $\ell,m\in (0,+\infty] $, $\vartheta\in \R $, the {\color{black}inhomogeneous} Besov space $B^\vartheta_{\ell,m}$  can be characterized {\color{black}(see \cite{Tri06})}  with
\[
{\color{black} B^\vartheta_{\ell,m}(\R^d)} = B^\vartheta_{\ell,m}=\left\{f\in\textcolor{black}{\mathcal S'}( \R^d)\,:\,\Vert f\Vert_{B^\vartheta_{\ell,m}}:=\|\mathcal F^{-1}(\phi\mathcal F(f))\|_{L^\ell}+\mathcal T_{\ell,m}^\vartheta(f)<\infty\right\},
\]
\begin{align}
	\mathcal T_{\ell,m}^\vartheta(f):=&\left\{
	\begin{aligned}
		&
		\left(\int_0^1\,\frac{\d v}{v}v^{(n-\vartheta/{\alpha})m}\Vert\partial^n_v p_\alpha(v,\cdot)\star f\Vert^m_{L^\ell}\right)^{\frac 1m}\,\text{ for }\,1\le m<\infty,\\
		&\sup_{v\in(0,1]}\left\{v^{(n-\vartheta/{\alpha})}\Vert\partial^n_v p_\alpha(v,\cdot)\star f\Vert_{L^\ell}\right\}\,\text{for}\,m=\infty,
	\end{aligned}
	\right. 
	\label{HEAT_CAR}
\end{align}
with $\star$ denoting the spatial convolution, $n$ being any non-negative integer (strictly) greater than $\vartheta/{\alpha}$, the function $\phi$ being a $\mathcal C^\infty_0$-function (infinitely differentiable function with compact support) such that $\phi(0)\neq 0$, and $p_\alpha(v,\cdot)$ denoting \textcolor{black}{as above} the density function at time $v$ of the $d$-dimensional isotropic stable process. 
{\color{black}The term $\mathcal T_{\ell,m}^\vartheta(f)$ appearing in the norm is called ``thermic part'' of the  Besov norm, and this characterization of the norm and related Besov spaces is \textcolor{black}{often} called the ``thermic characterization''.}

The controls needed for the analysis of the last time step of the error involving the Euler scheme (because of the frozen time in the dynamics) are the following:
\begin{itemize}
\item Duality inequality: for all $\vartheta \in \R$, $(\ell,m)\in [1,+\infty]^2$, with $m'$ and $\ell'$ respective conjugates of $m$ and $\ell$, and $(f,g)\in B_{\ell,m}^\vartheta \times B_{\ell',m'}^{-\vartheta}$,
	\begin{equation}\label{dual-ineq}
		\textcolor{black}{|\langle f,g\rangle_{B_{\ell,m}^{\vartheta},B_{\ell',m'}^{-\vartheta}}|:=}\left| \int f(y)g(y) \mathrm{d}y \right| \leq \Vert f \Vert_{B_{\ell,m}^\vartheta} \Vert g \Vert_{B_{\ell',m'}^{-\vartheta}},
	\end{equation}
	\textcolor{black}{where the duality pairing is here denoted in integral form for notational convenience}.
	We refer to Proposition 6.6 in \cite{lema:02} for a proof.

\item For $\vartheta\in (0,1) $, $\forall 0< s  < t$, $\forall (x,y)\in (\R^d)^2$, $\forall \zeta \in (\vartheta,1]$, $\forall k \in \{ 0,1\}$,
			\begin{align}\label{BESOV_DENSITY_LP}
			\| {p}_{\alpha} (s,x-\cdot) \nabla^k_y p_{\alpha} (t-s,y-\cdot) \|_{B^{\vartheta}_{p',q'}} \lesssim 
				\frac{\bar{p}_{\alpha} (t,x-y)}{(t-s)^{\frac{k}{\alpha}}} t^{-\frac{\vartheta}{\alpha}}\left[ \frac{1}{s^{\frac{d }{\alpha  p}}}+\frac{1}{(t-s)^{\frac{d }{\alpha  p}}} \right] \left[
				\frac{t^{\frac{\zeta}{\alpha}}}{s^{\frac{\zeta }{\alpha}}}+\frac{t^{\frac{\zeta}{\alpha}}}{(t-s)^{\frac{\zeta }{\alpha}}}  \right].
			\end{align}
			The proof of \eqref{BESOV_DENSITY_LP} can be found in \cite{Fit23} when $p_\alpha(s,x-\cdot) $ is replaced by $\bar p_\alpha(s,x-\cdot) $ in the equation. The proof can be extended to the considered case since the smoothness and controls required therein are shared by $p_\alpha $ and $\bar p_\alpha $.
\end{itemize}

\bibliographystyle{alpha}
\bibliography{ar}
\end{document}